\documentclass[reqno]{amsart}
\usepackage{a4wide}

\numberwithin{equation}{section}

\def\res{\mathop{\rm res}\limits}
\def\endproof{\bigskip}

\def\F{{\mathcal F}}
\def\L{{\mathcal L}}

\def\N{{\mathcal N}}

\def\C{{\mathcal C}}
\def\D{{\mathcal D}}

\newtheorem{theorem}{Theorem}[section]
\newtheorem{prop}[theorem]{Proposition}
\newtheorem{lemma}[theorem]{Lemma}
\newtheorem{corollary}[theorem]{Corollary}
\newtheorem{example}[theorem]{Example}

\newtheorem{Def}[theorem]{Definition}

\newdimen\tableauside\tableauside=1.0ex
\newdimen\tableaurule\tableaurule=0.4pt
\newdimen\tableaustep
\def\phantomhrule#1{\hbox{\vbox to0pt{\hrule height\tableaurule width#1\vss}}}
\def\phantomvrule#1{\vbox{\hbox to0pt{\vrule width\tableaurule height#1\hss}}}
\def\sqr{\vbox{%
  \phantomhrule\tableaustep
  \hbox{\phantomvrule\tableaustep\kern\tableaustep\phantomvrule\tableaustep}%
  \hbox{\vbox{\phantomhrule\tableauside}\kern-\tableaurule}}}
\def\squares#1{\hbox{\count0=#1\noindent\loop\sqr
  \advance\count0 by-1 \ifnum\count0>0\repeat}}
\def\tableau#1{\vcenter{\offinterlineskip
  \tableaustep=\tableauside\advance\tableaustep by-\tableaurule
  \kern\normallineskip\hbox
    {\kern\normallineskip\vbox
      {\gettableau#1 0 }%
     \kern\normallineskip\kern\tableaurule}%
  \kern\normallineskip\kern\tableaurule}}
\def\gettableau#1 {\ifnum#1=0\let\next=\null\else
  \squares{#1}\let\next=\gettableau\fi\next}
\begin{document}

\title[Frobenius manifolds: natural submanifolds \& induced bi-Hamiltonian structures]
{Frobenius manifolds: natural submanifolds and \\ induced bi-Hamiltonian structures}

\author{I.A.B. Strachan}
\date{\today}
\address{Department of Mathematics\\ University of Hull\\
Hull HU6 7RX\\ U.K.}

\email{i.a.strachan@hull.ac.uk}

\begin{abstract}
Submanifolds of Frobenius manifolds are studied. In particular, so-called natural submanifolds are defined
and, for semi-simple Frobenius manifolds, classified. These carry the structure of a Frobenius
algebra on each tangent space, but will, in general, be curved.
The induced curvature is studied, a main result being that these natural submanifolds
carry a induced pencil of compatible metrics. It is then shown how one may constrain the
bi-Hamiltonian hierarchies associated to a Frobenius manifold to live on these natural
submanifolds whilst retaining their, now non-local, bi-Hamiltonian structure.
\end{abstract}

\maketitle

\tableofcontents

\bigskip

\section{Introduction}

The study of the structures induced on a submanifold and their relationship
to the ambient manifold is one of the oldest problems in differential
geometry. The aim of this paper is to study the properties of submanifolds
of Frobenius manifolds.
Frobenius manifolds have a particularly rich structure, and have their
origin, as well as applications, in a wide range of seemingly disparate
areas of mathematics such as:
\begin{itemize}
\item{} topological quantum field theory;
\item{} algebraic/ennumerative geometry and quantum cohomology;
\item{} singularity theory;
\item{} integrable systems.
\end{itemize}
The emphasis in this paper will be on the purely geometric properties
of submanifolds and their application within the theory of integrable
systems. This will draw on ideas from singularity theory, in particular properties
of discriminants and caustics. Whether or not the ideas are
of relevance to the other areas is an open question.

\medskip

The key property of a Frobenius manifold is the existence of a Frobenius
algebra on each tangent space to the manifold:

\begin{Def}
An algebra $(\mathcal{A},\circ,<,>)$ over $\mathbb{C}$ is a Frobenius algebra if:
\begin{itemize}
\item{} the algebra $\{\mathcal{A},\circ\}$ is commutative, associative with unity $e\,;$
\item{} the multiplication is compatible with a $\mathbb{C}$-valued bilinear, symmetric, nondegenerate
inner product
\[
<,>\,: \,\mathcal{A}\times \mathcal{A}\rightarrow\mathbb{C}
\]
in the sense that
\[
<a \circ b,c> = < a,b\circ c>
\]
for all $a,b,c\in\mathcal{A}\,.$
\end{itemize}
\end{Def}
\noindent With this structure one may defined a Frobenius manifold \cite{D}:
\begin{Def} $(M,\circ,e,<,>,E)$ is a Frobenius manifold if each tangent
space $T_pM$ is a Frobenius algebra varying smoothly over $M$ with the
additional properties:
\begin{itemize}
\item{} the inner product is a flat metric on $M$ (the term \lq metric\rq~will denote a
complex-valued quadratic form on $M$);
\item{} $\nabla e=0$, where $\nabla$ is the Levi-Civita connection of the metric;
\item{} the tensor $(\nabla_W c)(X,Y,Z)$ is totally symmetric for all vectors
$W,X,Y,Z \in TM\,;$
\item{} a vector field $E$ must be determined such that
\[
\nabla(\nabla E) = 0\]
and that the corresponding one-parameter group of diffeomorphisms acts by
conformal transformations of the metric and by rescalings on the Frobenius
algebras $T_pM\,.$
\end{itemize}
\end{Def}
It is immediately apparent that an arbitrary submanifold $N\subset M$ of a Frobenius manifold
will not be a Frobenius manifold, as the induced metric will in general be curved.
Moreover, the induced multiplication on the subtangent space $T_pN\subset T_pM\,,p\in N$
will not be, in general, associative. Rather than develop a full structural theory
for submanifolds - which could easily be done - only so-called
natural submanifolds will be studied. On such submanifolds the induced multiplication
is associative and compatible with the induced metric. For semi-simple Frobenius
manifolds such natural submanifolds may be classified. The simplest example
comes from the Frobenius manifold constructed from the Coxeter group $A_3$. Here
the natural submanifolds are the swallow-tail discriminant, the cylinder over the
semi-cubical caustic and the planar Maxwell set.

The motivation for studying submanifolds came from two main
examples, more detail of which are given below. One of the best
understood classes of Frobenius manifolds come from the unfolding
of the $A_n$ singularity \cite{D}, $z^{n+1}\mapsto z^{n+1} + a_1 z^{n-1} +
\ldots a_n=p(z)\,.$ This derivation assumes that the roots of
$p^\prime(z)=0$ are distinct. However without this assumption,
i.e. with multiple roots, much of the structure of a Frobenius
manifold remains - one has a semi-simple Frobenius algebra on each
tangent space, compatible with an Euler vector field, and having a
covariantly constant identity vector field. However the metric is
no longer, in general, flat (a similar question was raised in
\cite{M}, III.7.1). This manifold of multiple roots should be thought
off a submanifold (in fact a caustic) of the original Frobenius
manifold. The details of this constitute Main Example A below. The
second motivation came from studying systems of hydrodynamic type
associated with Toda/Benney hierarchies \cite{FS,S1}, the simplest being
\begin{eqnarray*}
u_T & = & u v_X\,,\\
v_T & = & v u_X
\end{eqnarray*}
which is just the familiar dispersionless Toda equation (and hence
related to the quantum cohomology of $\mathbb{C}P^1$). An obvious
reduction of this system is to constrain the system to the
submanifold $u-v=0\,,$ reducing the system to the Riemann equation
$u_T = u u_X\,.$ Such submanifolds are clearly very special and
should be thought off a submanifold (in fact a discriminant) of the
original Frobenius manifold. The details of this constitute Main
Example B below.

\bigskip

\noindent{\bf Main Example A:} \cite{D} Consider the space $M$ of complex
polynomials
\[
p(z) = z^{m+1} + a_1 z^{m-1} + \ldots + a_m\,.
\]
Such a space carries the structure of a Frobenius manifold,
associated with the Coxeter group $A_m\,.$ Tangent vector to
$M$ take the form
\[
\dot p(z) = \dot a_1 z^{m-1} + \ldots + \dot a_m\,,
\]
and the algebra on the tangent space is
\[
A_p = \mathbf{C}[z] /p^\prime(z)
\]
and the inner product is
\begin{equation}
<f,g>_p = \res_{z=\infty} \left\{\frac{f(z) g(z)}{p^\prime(z)}\right\}\,.
\label{Anmetric}
\end{equation}
In terms of canonical coordinates $u^i=p(\alpha_i)$ where the
$\alpha_i$ are (distinct) roots of $p^\prime(z)=0$ the metric
becomes diagonal and Egoroff:
\begin{eqnarray}
g & = & \sum \frac{1}{p^{\prime\prime}(\alpha_i)} (du^i)^2\,,\label{simplepoles}\\
\frac{1}{p^{\prime\prime}(\alpha_i)} & = & \frac{1}{m+1}
\frac{\partial a_1}{\partial u^i}\,.
\label{simplepoles2}
\end{eqnarray}
Note that this all assumes that the roots of $p^\prime(z)=0$ are
distinct. In what follows no such assumption will be assumed. The
generic features will remain - the metric will remain diagonal and
Egoroff. The result will be incorporated in a more general scheme
which will be developed over the subsequent sections. One may
regard the manifold $N$ with repeated roots as a submanifold (in fact
a caustic) in $M\,.$

Suppose that
\[
p^\prime(z) = (m+1) \prod_{i=1}^n (z-\alpha_i)^{k_i}
\]
where $k_i\geq 1\,, \sum^n_{i=1} k_i=m\,, \sum^n_{i=1} k_i \alpha_i=0\,.$
Coordinates on $N$ are defined by $\tau^i=p(\alpha_i)\,.$ It
follows immediately that
\begin{eqnarray*}
\delta_{ij} & = & \left. \frac{\partial p}{\partial\tau^j}
\right|_{z=\alpha_i} \,, \\
0 & = & \left.\frac{d^k~}{d z^k}
\frac{\partial p}{\partial\tau^j} \right|_{z=\alpha_i}\,,
\quad\quad 1 \leq k \leq k_i-1\,.
\end{eqnarray*}
A simple parameter counts yields $m$ equations for the $m$
unknowns in $\frac{\partial p}{\partial\tau^j}\,.$ This gives two
different forms for $\frac{\partial p}{\partial\tau^j}\,:$
\begin{eqnarray}
\frac{\partial p}{\partial\tau^j} & = & p^{(j)}(z) \prod_{i\neq j}
(z-\alpha_i)^{k_i}\,, \label{dpdtau1}\\
& = & 1 + \sum_{k=k_j}^{m-1} \frac{1}{k!} \left.
\frac{d^k~}{dz^k} \left( \frac{\partial p}{\partial \tau^j}\right)
\right|_{z=\alpha_j} \,. (z-\alpha_j)^k\,,\label{dpdtau2}
\end{eqnarray}
where $p_j$ is a polynomial of degree $k_j-1$ with
$p_j(\alpha_j)\neq 0\,.$ The metric on $N$ is given by
\[
g_{ij} = \res_{z=\infty} \left\{ \frac{ \frac{\partial p}{\partial\tau^i}
\frac{\partial p}{\partial\tau^j}}{p^\prime}\right\} dz
\]
and using (\ref{dpdtau1}) gives immediately that the metric is
diagonal and, on using (\ref{dpdtau2}), that
\[
g_{ii} = \res_{z=\alpha_i} \left\{\frac{1}{p^\prime(z)}\right\}\,.
\]
In the case of simple poles this reduces to (\ref{simplepoles}).
An immediate corollary of this is that
\begin{equation}
\sum_{i=1}^n g_{ii} = \frac{1}{2\pi i} \oint_{\mathcal C}
\frac{dz}{p^\prime(z)}=0\,, \label{summetric}
\end{equation}
where $\mathcal C$ is a large contour containing all of the
$\alpha_i\,.$ To show that this metric is Egoroff is considerable
more involved, even though the final result is simple. A more
geometric proof will be given below, here it will be derived by
direct calculation.

Let
\begin{eqnarray*}
h^{(i)} & = & \prod_{r\neq i} (z - \alpha_r)^{k_r} \,,\\
& = & \sum_{s=0}^{m-k_i} h_s^{(i)} \frac{ (z-\alpha_i)^s}{s!}
\end{eqnarray*}
and define the coefficients $h^{(i,-1)}_s$ by the inverse series
\[
(h^{(i)})^{-1} =
\sum_{s=0}^{\infty} h_s^{(i,-1)} \frac{ (z-\alpha_i)^s}{s!}\,.
\]
Since $h_0^{(i)} = \prod_{r\neq i} (\alpha_i - \alpha_r)^{k_r}
\neq 0$ these coefficients are uniquely defined. To proceed
further
one requires the following:

\begin{lemma} Consider the expansion of $p^{(i)}(z)$ around
$z=\alpha_i\,,$ and let
\[
p_r^{(i)} = \left.\frac{d^r p^{(i)}}{dz^r}\right|_{z=\alpha_i}
\]
then $p_r^{(i)} = h^{(i,-1)}_r$ for $r=0\,, \ldots\,, k_i-1\,.$
\end{lemma}

\noindent This is proved by showing that the linear equations for
the $h^{(i,-1)}_r$ and the $p_r^{(i)}$ are identical. In the case
of simple zeros this result is immediate. Note that the explicit
form of these coefficients is not required, just their equality.
Then
\begin{eqnarray*}
g_{ii} & = & \res_{z=\alpha_i} \frac{1}{p^\prime(z)} \,, \\
& = & \frac{1}{m+1} \res_{z=\alpha_i} \frac{1}{(z-\alpha_i)^{k_i}}
(h^{(i)})^{-1}\,, \\
& = &
\frac{1}{m+1} \res_{z=\alpha_i} \frac{1}{(z-\alpha_i)^{k_i}}
\sum_{s=0}^{\infty} h_s^{(i,-1)} \frac{ (z-\alpha_i)^s}{s!}\,,\\
& = & \frac{1}{m+1} \, \frac{h^{(i,-1)}_{k_i-1}}{(k_i-1)!}\,,\\
& = & \frac{1}{m+1} \, \frac{p^{(i)}_{k_i-1}}{(k_i-1)!}\,,\\
& = & \frac{1}{m+1}
{\rm~coefficient~of~}z^{m-1}{\rm~in~expansion~of~}
\frac{\partial p}{\partial\tau^i}\,, \\
& = & \frac{1}{m+1} \frac{\partial a_1}{\partial\tau^i}\,.
\end{eqnarray*}
Hence the metric is Egoroff. Other properties may be similarly
derived; $N$ carries an Euler vector field and a covariantly
constant unity vector field $e$ (this following from
(\ref{summetric}))\,.

\bigskip

\endproof

\noindent{\bf Main Example B:} The multicomponent Toda hierarchy is defined in terms of a Lax function

\[
\L(z) = z^{M-1} + \sum_{i=-1}^{M-2} z^i S^i(X,{\bf T})\,,
\quad\quad{\bf T} =\{T_1\,,T_2\,,\dots \}
\]

\noindent by the Lax equation

\begin{equation}
\frac{\partial\L}{\partial T_n} = \{ (\L^{\frac{n}{M-1}})_{+} , \L \}\,.
\label{eq:lax}
\end{equation}

\noindent Here the bracket is
defined by the formula

\[
{ \{ f,g \} } = z \frac{\partial f}{\partial z}
\frac{\partial g}{\partial X} - z \frac{\partial f}{\partial X}
\frac{\partial g}{\partial z}\,.
\]

\noindent and $( {\mathcal O} )_{+}$ denotes the projection of the
function $\mathcal O$ onto non-negative powers of $p\,.$ For
example, one obtains from the Lax equation (\ref{eq:lax}) with
$M=2\,,n=1$ the system (where $\L=z+S+P z^{-1})\,$
\begin{equation}
\begin{array}{rcl}
S_T & = & P_X \,, \\
P_T & = & P S_X \,.
\end{array}
\label{eq:toda}
\end{equation}
and, with $M=3\,,n=1$ the system (where $\L=z^2 + Sz + P + Q z^{-1}\,$)
\begin{equation}
\begin{array}{rcl}
S_T & = & P_X - \frac{1}{2} S S_X \,, \\
P_T & = & Q_X \,, \\
Q_T & = & \frac{1}{2} Q S_X\,.
\end{array}
\label{eq:3toda}
\end{equation}
A change of dependent variables from the $\{ S^i(X,{\bf T}) \}$ to
so-called modified variables $\{v^i(X,{\bf T}) \}$ defined by a
factorization of the Lax equation
\[
\L=\frac{1}{z} \prod_{i=1}^N [z+v_i(X,{\bf T}) ]
\]
provides an extremely useful computational tool in the study of
the Toda hierarchy\cite{FS,S1}. Quantities such as
\[
Q^{(n)} = \frac{1}{2\pi i} \oint \L^{ \frac{n}{
m-1} } \frac{dz}{z}
\]
which are are conserved with respect to the evolutions defined by the
Lax equation (\ref{eq:lax}) may be evaluated for all values of
$M$ and $n\,$ in terms of a simple combinatorial formula
\[
Q^{(n)} = \sum_{ \{ r_i \,: \sum_{i=1}^M r_i= n \} } \Bigg\{
\prod_{i=1}^M
\begin{pmatrix}
\frac{n}{M-1} \cr r_i
\end{pmatrix}
v_i^{r_i} \Bigg\}\,,
\]
and similar formulae exist for the evolution equations themselves.
Thus the modified variables enables one to perform the
general calculations with an arbitrary numbers of fields with
little increase in complexity. The geometrical significance
of these variables is that they are
basically the flat coordinates for the intersection form of the
underlying Frobenius manifold, or equivalently, flat coordinates
for the second Hamiltonian structure. This Hamiltonian structure is
defined by
the manifestly flat metric
\[
{\bf g} = \sum_{i \neq j} \frac{dv_i}{v_i} \,
\frac{dv_j}{v_j}\,
\]
(so the actual flat coordinates are ${\tilde v}_i = \log v_i\,).$

\bigskip

Consider (\ref{eq:toda})
written in terms of these modified variables
\begin{eqnarray*}
u_T & = & u v_X\,,\\
v_T & = & v u_X\,
\end{eqnarray*}
where $S=u+v$ and $P=uv\,.$
Similarly the 3-component system (\ref{eq:3toda}) transforms to
\begin{eqnarray*}
u_T & = &  u ( -u_X+v_X+w_X ) \,,\\
v_T & = &  v ( +u_X-v_X+w_X ) \,,\\
w_T & = &  w ( +u_X+v_X-w_X ) \,.
\end{eqnarray*}
where $S=u+v+w\,,P=uv+vw+wu$ and $Q=uvw\,.$ It is clear from the symmetric form of these
equations that one possible reduction is to constrain the systems onto the surface given by
the constraint $u-v=0\,.$ In terms of the original variables this
corresponds to the constraint $S^2-4P=0$ and in the second this
corresponds to the constraint
\[
4P^3+27 Q^2 - 18 PQS - P^2 S^2 + 4 Q S^3 = 0 \,,
\]
these being the condition for the corresponding polynomial equation
$\L(z)=0$ to have a double root. Clearly these ideas generalize to an arbitrary
number of fields and arbitrary multiple roots. Geometrically one
is constraining an $M$ dimensional system onto a $N$ dimensional
submanifold. Note that these reductions are far
easier to study using these modified variables. All these results
generalize to rational Lax equations in an
entirely analogous fashion.

\bigskip

\endproof

These two examples provided the motivation for the study of submanifolds. The
key property possessed by the submanifolds in both these examples is the
commutative, associative and quasihomogeneous multiplication on the
subspace's tangent bundle. Submanifolds with such induced structures will be referred to
as \lq natural\rq~submanifolds. These examples also have an induced identity vector
field on the submanifolds, and hence one has a Frobenius algebra on each
tangent space of the submanifold. This paper is a much extended version of the paper
\cite{S3}.

\bigskip

\section{$F$-manifolds and their natural submanifolds}

The definition of a Frobenius manifold consists of a large number
of intermeshing parts, and it is perhaps difficult to see which components
are the most important. One point of view, coming via singularity theory,
is that it is the multiplication and the Euler vector field which are the
central elements; the existence of a compatible flat metric being, for example,
derived results. This point of view is encapsulated in the weaker notion of
an $F$-manifold \cite{HM,He,M}. Here one has a commutative, associative multiplication
with a single additional property, automatically satisfied in the case of Frobenius manifolds.
Starting with an $F$-manifold one may gradually add additional structures and compatibility
conditions until one obtains a Frobenius manifold. This has the advantage that one
may see on what structures the various compatibility conditions depend.

A similar approach will be taken here for submanifolds. One may defined
a natural submanifold $N$ of an $F$-manifold $M$ by requiring that $TN\circ TN\subset TN\,.$
As one adds various structures and compatibility conditions onto $M$ one can
also study the induced structures
on $N$ and the failure, of otherwise, of the associated compatibility conditions.
The various results in this section are formulated with this approach in mind, even though
the main aim is to study natural submanifolds of Frobenius manifolds. The results in
this section are mainly algebraic; curvature properties being studied in section 3.

\bigskip

\begin{Def} {\cite{HM,M}} An $F$-manifold is a pair $(M,\circ)$
where $M$ is a manifold and $\circ$ is a commutative, associative
multiplication $\circ\,:\,TM \times TM\rightarrow TM$ satisfying
the following condition:

\begin{equation}
Lie_{X\circ Y} (\circ) = X \circ Lie_Y(\circ) + Y \circ Lie_X
(\circ)\,,\quad\quad \forall X,Y\,\in TM\,.
\label{Fmandef}
\end{equation}
\end{Def}
\noindent Expanding the definition yields the equivalent condition
\begin{eqnarray*}
[X\circ Y, Z\circ W] - [X\circ Y,Z]\circ W - [X\circ Y,W]\circ Z&&\\
-X\circ [Y,Z\circ W] + X\circ [Y,Z]\circ W + X\circ [Y,W]\circ Z&&\\
-Y\circ [X,Z\circ W] + Y\circ [X,Z]\circ W + Y\circ [X,W]\circ Z&=&0
\end{eqnarray*}
for all $W,X,Y,Z\in TM\,.$
To such a manifold one may add various structures, demanding that
they are compatible with the multiplication.

\begin{Def} (a) An $F_E$ manifold is an $F$-manifold with an Euler
field of weight $d\,.$ This is a global vector field satisfying
the condition
\begin{equation}
Lie_E(\circ) = d \cdot\circ\,.
\label{Edef}
\end{equation}
\medskip
\noindent (b) An $F_g$ manifold is an $F$-manifold with a metric
$<,>$ compatible with the multiplication:
\begin{equation}
<X\circ Y,Z> = <X,Y\circ Z> \,, \quad\quad X,Y,Z\in TM\,.
\label{metricdef}
\end{equation}
\medskip
(c) An $\F$ manifold is both an $F_E$ and an $F_g$ manifold, with
the $E$ and $g$ related by the relation
\begin{equation}
Lie_E \, <,> = D \, <,>\,
\label{Fmandef2}
\end{equation}
for some constant $D\,.$
\end{Def}
Expanding definition (\ref{Edef}) yields the equivalent condition
\[
[E,X\circ Y] - [E,X] \circ Y - X\circ[E,Y] - d\cdot X\circ Y = 0
\]
for all $X,Y\in TM\,,$ and (\ref{Fmandef2}) yields the equivalent condition
\[
E<X,Y> - <[E,X],Y> - <X,[E,Y]> = D <X,Y>
\]
for all $X,Y\in TM\,.$

The following definition of a natural submanifold will play a
central role in this paper.

\begin{Def} A natural submanifold $N$ of an $F_E$ manifold
$(M,\circ,E)$ is a submanifold $N\subset M$ such that:
\begin{itemize}
\item[(a)] $TN \circ TN \subset TN\,,$
\item[(b)] $E_x \in TN$ for all $x\in N$\,.
\end{itemize}
\end{Def}
\noindent One could clearly define the notion of a natural submanifold of an
$F$-manifold by ignoring the second condition. An immediate
consequence of this definition is the following basic result, the
proof of which follows from the fact that if $X\,,Y\in TN$ then $[X,Y]\in TN$:

\begin{lemma} All natural submanifolds of an $F_E$ manifold are
$F_E$ manifolds with respect to the natural induced structures.
\end{lemma}

\medskip

\begin{example} {\rm (Massive $F$-manifolds) Given a semi-simple
$F_E$ manifold $(M,\circ,E)$ the tangent space $T_pM$ at a generic point
decomposes into into one-dimensional algebras with
\[
\delta_i \circ \delta_j = \delta_{ij} \delta_i
\]
(so $\delta_i$ are the idempotents of the algebra on $T_pM$). The
$F$ manifold condition ensures that these vector fields
commute $[\delta_i,\delta_j]=0$ and hence provide a
canonical coordinate system $\{u^i \}$ with
\[
\partial_i =\frac{\partial~}{\partial u^i}\,.
\]
In this basis one has then:
\begin{eqnarray*}
\frac{\partial~}{\partial u^i} \circ \frac{\partial~}{\partial u^j}
& = & \delta_{ij} \frac{\partial~}{\partial u^i}\,, \quad\quad
i,j=1\,,\ldots\,,m={\rm dim}\,M\,,\\
E & = & \sum_{i=1}^{{\rm dim}\,M}  u^i\frac{\partial~}{\partial
u^i}\,.
\end{eqnarray*}
Then the submanifolds defined by the level sets
\[
\underbrace{\{\quad u^i = 0 \,, i \in \D\quad\}}_{\rm discriminant~hypersurfaces} \cap\,
\underbrace{\{\quad u^i - u^j = 0 \,, (i,j) \in \C \quad\}}_{\rm
caustic~hypersurfaces}
\]
are natural $F_E$ manifolds. Here $\D$ and $\C$ are arbitrary subsets of $I$ and $I\times I\,$
where $I=\{1\,,\ldots\,,{\rm dim}\,M\}\,.$ This example will turn out to be
canonical for natural submanifolds  semi-simple Frobenius
manifolds.

Relabelling the coordinates gives the following parametrization of
a natural submanifold:
\[
(u^1\,,\ldots\,,u^m)=(
\underbrace{\tau^1\,,\ldots\,,\tau^1}_{k_1};\ldots;
\underbrace{\tau^n\,,\ldots\,,\tau^n}_{k_n};
\underbrace{0\,,\ldots\,,0}_{m-(k_1+\ldots + k_n)})\,.
\]
In what follows the notation $(k_1\,,\ldots\,,k_n\,,0)$ will be
used to denote a particular submanifold, so the original
manifold would just
be $(1\,,1\,\ldots\,,1)\,.$
An alternative notation is to use a Young tableau, so, for example,
the $(4,2)$ caustic would be denoted $\tableau{4 2}\,.$
The terms pure discriminant will refer
to a submanifold where $k_i=1$ and a pure caustic will refer to a
submanifold where $\sum k_i = m\,.$
For ${\dim}\,M=2$ the only possibilities
are (the notation $N \longrightarrow M$ means that
$N$ is a natural codimension one submanifold of $M$):

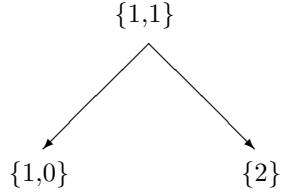
\begin{picture}(100,100)(-120,0)
\put(40,80) {\{1,1\}}
\put(53,72) {\vector(1,-1){40}}
\put(53,72) {\vector(-1,-1){40}}
\put(0,20) {\{1,0\}}
\put(88,20) {\{2\}}
\end{picture}

\noindent For ${\dim}\,M=3$ one obtains the following strata of nested
submanifolds:

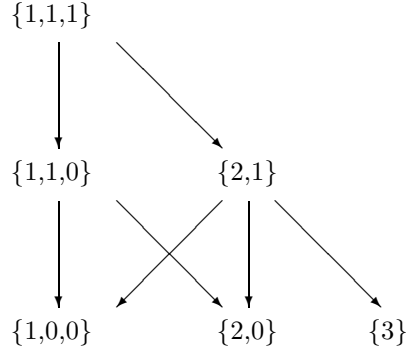
\begin{picture}(250,180)(-50,0)
\put(50,140){\{1,1,1\}} \put(68,132) {\vector(0,-4){40}}
\put(50,80){\{1,1,0\}} \put(68,72) {\vector(0,-4){40}}
\put(50,20){\{1,0,0\}}
\put(90,132) {\vector(1,-1){40}}
\put(90,72) {\vector(1,-1){40}}
\put(130,72) {\vector(-1,-1){40}}
\put(128,80) {\{2,1\}}
\put(128,20) {\{2,0\}}
\put(140,72){\vector(0,-4){40}}
\put(150,72){\vector(1,-1){40}}
\put(185,20) {\{3\}}
\end{picture}

\noindent For ${\dim}\,M>3$ such diagrams become considerably more
complicated, the number of such submanifolds being $\sum_{n=1}^m
\mu(n)\,,$ where $\mu(n)$ is the number of partitions of $n\,.$}
\end{example}

\endproof

\bigskip

Suppose now one has an $F_g$-manifold. Then on any (non-null)
submanifold $N$ one may define an induced metric $g_N$ and
also an induced product $\star\,: TN \times TN \rightarrow TN$
where $\star$ is defined by
\[
X\star Y = pr(X \circ Y) \quad \forall X\,,Y\,\in T_x N \subset T_x M\,,
\]
where $pr$ denotes the projection (using the original metric $g$ on $M$) of
$u\circ v\in T_x M$ onto $T_x N\,.$ This induced multiplication may have very
different algebraic properties than those of its progenitor.
However the induced metric and multiplication remain compatible.

\begin{lemma}
The induced structures satisfy the condition
\[
<X \star Y,Z> = <X, Y \star Z> \quad\quad \forall X,Y,Z,\in T_x N\,.
\]
\end{lemma}

\noindent The proof following immediately from the definitions.
Putting these results together gives the following:

\begin{prop} Any natural submanifold of an $\F$-manifold is an
$\F$-manifold with respect to the naturally induced structures.
\end{prop}

\noindent One note of caution though: this is a formal result - it may be
the case that the induced metric is not defined or is degenerate
on a specific natural submanifold.

\medskip

The definition of a natural submanifold just uses the multiplication
and the Euler vector field. If one has, in addition, an identity
vector field then a natural submanifold will inherit an induced
identity field:

\begin{lemma} Let $(M,\circ,g)$ be an $F_g$ manifold with a unity vector
field $e$ and let $N$ be a natural submanifold of $M$. Then $N$ possesses
an induced identity vector field.
\end{lemma}

\medskip

\noindent{\bf Proof~} (Note, this lemma only uses part (a) of the definition of
a natural submanifold.) Using the metric, one has an orthogonal decomposition
(assuming the induced metric on $N$ is not degenerate) of the tangent space
$T_xM$ at points $x\in N\,$):
\begin{equation}
T_x M  \cong (T_x N) \,\,\oplus \,\,(T_x N)^\perp\,,\quad \quad x\in N\,,
\label{tangentspacedecomp}
\end{equation}
so $e$ decomposes as $e=e^\top+e^\perp\,.$ Hence
\[
X\circ e^\perp = X-X \circ e^\top \in T_xN\,.
\]
Clearly $<X\circ e^\perp,n>=0$ for all $n\in (T_xN)^\perp$ and
\begin{eqnarray*}
<X\circ e^\perp,Y> & = & <X \circ Y, e^\perp> \,,\\
& = & 0
\end{eqnarray*}
for all $Y\in T_xN\,,$ using the invariance property of the multiplication.
Thus $X\circ e^\perp = 0 $ and hence $X\circ e^\top = X\,$ for all $X \in T_xN\,.$

\endproof

\noindent An immediate corollary of this is:

\begin{corollary} Let $M$ be a Frobenius manifold and let $N\subset M$ be a natural
submanifold. Then each tangent space $T_xN$ carries the structure of a
Frobenius algebra with respect to the induced structures.
\end{corollary}

\bigskip

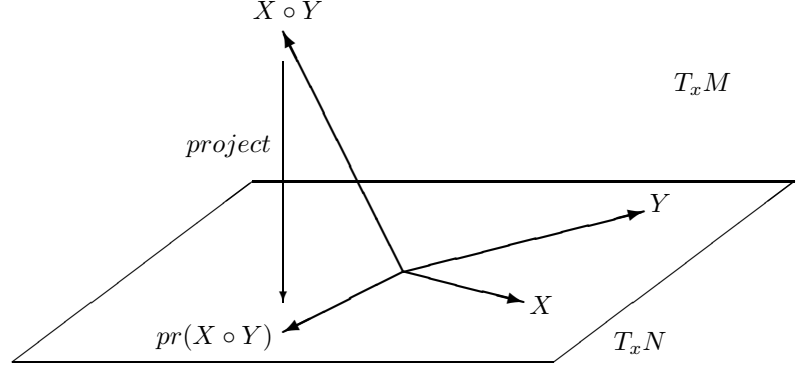
\begin{figure}

\setlength{\unitlength}{0.8cm}
\begin{picture}(17,7)
\put(3,0){\line(1,0){9}}
\put(12,0){\line(4,3){4}}
\put(3,0){\line(4,3){4}}
\put(7,3){\line(1,0){9}}
\thicklines
\put(9.5,1.5){\vector(4,-1){2}}
\put(9.5,1.5){\vector(4,1){4}}
\put(9.5,1.5){\vector(-1,2){2}}
\put(9.5,1.5){\vector(-2,-1){2}}
\put(11.6,0.8){$X$}
\put(13.6,2.5){$Y$}
\put(7.0,5.7){$X\circ Y$}
\thinlines
\put(7.5,5.0){\vector(0,-1){4}}
\put(5.4,0.3){$pr(X\circ Y)$}
\put(5.9,3.5){$project$}
\put(13,0.2){$T_x{N}$}
\put(14,4.5){$T_x{M}$}
\end{picture}

\caption{The definition of the induced multiplication}
\end{figure}

For a semi-simple $\F$-manifold one may classify all natural
submanifolds, at least formally. The idea is to describe an
arbitrary submanifold as the intersection of level sets,
$N=\cap \{ \phi^{{\tilde\alpha}} = 0 \}\,,$ the geometric
conditions on $N$ to be a natural submanifold then reduce to
a simple set of overdetermined
partial differential equations for the functions
$\phi^{{\tilde\alpha}}$
which may be solved.

\begin{theorem} Let $\{M,\circ,E,g\}$ be a semi-simple
$\F$ manifold. Then:
\begin{itemize}
\item[{(a)}] the only natural submanifolds are those given
in the above example;
\item[{(b)}] the identity field is tangential to a natural
submanifold if and only if it is a pure caustic.
\end{itemize}
\end{theorem}

\noindent{\bf Proof } (a) Let $\imath : N\rightarrow M$ be the
inclusion of a submanifold $N$ in the manifold $M\,.$ Vector
fields on $N$ may be pushed-forward to vector fields on $M\,.$
Adopting a parametrization of the submanifold $N\,,$ so $u^i =
u^i(\tau^\alpha)\,,$ where
$i=1\,,\ldots\,,m\,,\alpha=1\,,\ldots\,,n\,,$ one obtains
\begin{eqnarray*}
\imath_\star : TN & \rightarrow & TM \,, \\
\imath_\star \left(\frac{\partial~}{\partial \tau^\alpha}\right) & = &
\frac{\partial u^i}{\partial \tau^\alpha}\frac{\partial~}{\partial u^i}\,.
\end{eqnarray*}
Similarly \cite{S2}, using the orthogonal decomposition (\ref{tangentspacedecomp})
(assuming the induced metric on $N$ is not degenerate):
\begin{equation}
\frac{\partial~}{\partial u^i} = A^\alpha_i
\frac{\partial~}{\partial \tau^\alpha} + n^{\tilde\alpha}_i
\frac{\partial~}{\partial \nu^{\tilde\alpha}} \label{orthogdecomp}
\end{equation}
where ${\rm span}(\partial_\nu) = (T_x N)^\perp\,.$

Consider now
\begin{eqnarray*}
\frac{\partial~}{\partial \tau^\alpha} \circ \frac{\partial~}{\partial \tau^\beta}
& = &
\frac{\partial u^i}{\partial \tau^\alpha}
\frac{\partial u^j}{\partial \tau^\beta}\frac{\partial~}{\partial u^i}\circ
\frac{\partial~}{\partial u^j}\,,\\
& = & \sum_i
\frac{\partial u^i}{\partial \tau^\alpha}
\frac{\partial u^i}{\partial \tau^\beta}
\left(
A^\alpha_i
\frac{\partial~}{\partial \tau^\alpha} + n^{\tilde\alpha}_i
\frac{\partial~}{\partial \nu^{\tilde\alpha}}
\right)\,,
\end{eqnarray*}
on using the canonical multiplication. To ensure that
$TN\circ TN \subset TN$ one must have

\begin{equation}
\Xi_{\alpha\beta}^{\tilde\alpha}=0
\label{obstruction}
\end{equation}

\noindent where
\[
\Xi_{\alpha\beta}^{\tilde\alpha}=
\sum_i
\frac{\partial u^i}{\partial \tau^\alpha}
\frac{\partial u^i}{\partial \tau^\beta}
n^{\tilde\alpha}_i
\,.
\]
To proceed further one adopts a Monge parametrization of $N$ so
\begin{eqnarray*}
u^i & = & \tau^i \,, \qquad\qquad i=1\,, \ldots \,, n\,, \\
u^{n+{\tilde\alpha}} & = & h^{\tilde\alpha} ( \tau^\alpha ) \,,
\quad\quad {\tilde\alpha} = 1\,,\ldots \,, m-n\,.
\end{eqnarray*}
With this $N$ may be described as the intersection of level sets
\[
N = \bigcap_{{\tilde\alpha}=1}^{m-n} \{ \phi^{\tilde\alpha} = 0 \}
\]
where $\phi^{\tilde\alpha} = h^{\tilde\alpha} -
u^{n+{\tilde\alpha}}\,.$
This may be used to find the normal vectors $n_i^{\tilde\alpha}$
with which the condition $\Xi_{\alpha\beta}^{\tilde\alpha}=0$
become
\[
\delta_{\alpha\beta}
\frac{\partial h^{\tilde\alpha}}{\partial\tau^\alpha} =
\frac{\partial h^{\tilde\alpha}}{\partial\tau^\alpha}
\frac{\partial h^{\tilde\alpha}}{\partial\tau^\beta}\,,
\quad\quad\alpha\,,\beta=1\,,\ldots\,,n\,,
{\tilde\alpha}=1\,, \ldots\,, m-n\,.
\]
If $\alpha=\beta$ then $h^{\tilde\alpha}_\alpha=0$ or $1\,.$ But
if $\alpha\neq\beta$ then
$h^{\tilde\alpha}_\alpha h^{\tilde\alpha}_\beta=0$ which
implies that $h^{\tilde\alpha}_\alpha=0$ except, possibly,
for one values of $\alpha\in\{1\,,\ldots\,,n\}\,.$ Such a value
will be denoted $\pi({\tilde\alpha})\,.$ Hence there are two
possibilities:
\[
h^{\tilde\alpha} = a^{\tilde\alpha} \,,\quad\quad\quad
h^{\tilde\alpha} = u^{\pi({\tilde\alpha})} + b^{\tilde\alpha}
\]
for arbitrary constants $a^{\tilde\alpha}\,,b^{\tilde\alpha}\,.$
Note, if one was to consider semi-simple $F_g$ manifolds, then the
classification of natural submanifolds would stop here.

For $N$ to be a natural submanifold requires the further condition
$E_x\in T_x N\,,$
\begin{eqnarray*}
E_x & = & \sum_{i=1}^m u^i \frac{\partial~}{\partial u^i}\,,\\
& = & \sum_{i=1}^m u^i\left(
A^\alpha_i
\frac{\partial~}{\partial \tau^\alpha} + n^{\tilde\alpha}_i
\frac{\partial~}{\partial \nu^{\tilde\alpha}}
\right)\,.
\end{eqnarray*}
Thus $(E_x)^\perp=0$ implies, using this parametrization, that
\[
\sum_{i=1}^m u^i \frac{\partial~}{\partial u^i} h^{\tilde\alpha} =
h^{\tilde\alpha} \,,
\]
so the $h^{\tilde\alpha}$ must be homogeneous functions of degree
$1\,.$ Hence $a^{\tilde\alpha}=b^{\tilde\alpha}=0\,.$ Thus
\begin{equation}
h^{\tilde\alpha} = 0 \,,\quad\quad\quad
h^{\tilde\alpha} = u^{\pi({\tilde\alpha})} \,.
\label{defh}
\end{equation}
On renaming the
coordinates one arrives at the examples described above.

\medskip

\noindent (b) Note that, as a consequence of semi-simplicity, there exists
a unity vector field
\[
e=\sum_{i=1}^m \frac{\partial~}{\partial u^i}
\]
with the property that $e \circ X=X$ for all $X\in TM\,.$ Similarly
\[
e_N = \sum_{\alpha=1}^n \frac{\partial~}{\partial \tau^\alpha}
\]
will be a unity vector field on $TN\,.$ Consider now the restriction of $e$ to
a natural submanifold. Using the above formulae, and in particular
(\ref{defh}), it is straightforward to show that
\[
\sum_{i=1}^m \frac{\partial~}{\partial u^i} =
\sum_{\alpha=1}^n \frac{\partial~}{\partial \tau^\alpha}-
\sum_{{\tilde\alpha}=1}^{m-n}
\left\{
1-\sum_{j=1}^n \frac{\partial h^{\tilde\alpha}}{\partial
\tau^i}\right\}
\frac{\partial~}{\partial \nu^{\tilde\alpha}}\,.
\]
Hence
$e^\bot=0$ if and only if
\[
\sum_{j=1}^n \frac{\partial h^{\tilde\alpha}}{\partial \tau^i}=1
\qquad\qquad\forall {\tilde\alpha}=1\,,
\ldots\,, m-n\,,
\]
that is, using (\ref{defh}), if and
only if $N$ is a pure caustic.

\endproof

\noindent In the next section curvature properties of natural submanifolds
will be examined.

\bigskip

\section{Frobenius manifolds and the curvature properties of
natural submanifolds}

Given an $\F$-manifold one may define the following tensors,
$c(X,Y,Z)=g(X\circ Y,Z)$, which, from (\ref{metricdef}), is a
totally symmetric $(3,0)$ tensor, $\nabla\circ$ and $\nabla c\,.$ The following
theorem is due to Hertling \cite{He}:

\begin{theorem} Let $(M,\circ,\nabla)$ be a manifold $M$ with a
commutative associative multiplication $\circ$ on $TM$ and with a
torsion free connection $\nabla\,.$ By definition,
$\nabla\circ(X,Y,Z)$ is symmetric in $Y$ and $Z\,.$ If the
$(3,1)$-tensor $\nabla\circ$ is symmetric in all three arguments,
then the multiplication satisfies for any local vector fields $X$
and $Y$
\[
Lie_{X\circ Y}(\circ) = X \circ Lie_Y (\circ) + Y \circ
Lie_X(\circ)\,.
\]
\end{theorem}
The converse, however, is false; one requires the properties of a
unity vector field. So far little mention has been made of the possibility of having a
unity vector field $e$ on $M\,,$ i.e. $e\in TM$ such that
$e\circ X=X\,,\forall X\in TM\,.$ Such fields play an important
role in Frobenius and $\F$ manifolds, as it connects the metric and the
multiplication since $<X,Y>=c(X,Y,e)\,.$ With this field one may
prove the following, again due to Hertling \cite{He}:

\begin{theorem} Let $(M,\circ,e,g)$ be a manifold with a
commutative and associative multiplication $\circ$ on $TM$, a unit
field $e$, and a metric $<,>$ on $TM$ which is multiplication
invariant (\ref{metricdef}). $\nabla$ denotes the Levi-Civita
connection fo the metric. The coidentity $\varepsilon$ is the
1-form defined by $\varepsilon(X)=<X,e>\,.$ The following
conditions are equivalent:
\begin{itemize}
\item[{(i)}] $(M,\circ,e)$ is an $F$ manifold and $\varepsilon$ is
closed;
\item[(ii)] the $(4,0)$ tensor $\nabla c$ is totally symmetric;
\item[(iii)] the $(3,1)$ tensor $\nabla \circ$ is totally
symmetric.
\end{itemize}
\end{theorem}

The property that $\nabla c$ is a totally symmetric
$(4,0)$ tensor is sometimes referred to as
quasi-potentially conditions since, if the metric is flat,
one may integrate the equations to give $c$ in terms of
derivatives of a prepotential $F\,,$ i.e.
\[
c(X,Y,Z) = XYZ (F)\,.
\]

\begin{Def} A Frobenius manifold is an $\F$ manifold
$(M,\circ,E,g)$ endowed with a unity vector field $e$
and satisfying the following conditions:
\begin{itemize}
\item[(a)] $g$ is flat;
\item[(b)] $d \varepsilon=0$ where
$\varepsilon(\cdot)=<e,\cdot>\,;$
\item[(c)] $Lie_e <,> =0\,;$
\item[(d)] $\Delta(U,V) = 0\,.$
\end{itemize}
where
$\Delta(U,V) = \nabla_U \nabla_V E - \nabla_{\nabla_U V} E\,,$
and $\nabla$ is the Levi-Civita connection.
\end{Def}
This definition differs somewhat from the conventional one given above.
Conditions (b) and (c) together imply that $\nabla e=0\,$ By the
above theorem, condition (b) implies that $\nabla c$ is totally
symmetric, and condition (a) then implies that there exists a
prepotential. Condition (d) perhaps needs a little explanation. If
$U$ and $V$ are flat vector fields the $\Delta(U,V) =\nabla_U
\nabla_V E \,.$ Alternatively, $\Delta_{ij}^k =\nabla_i \nabla_j
E^k$ in terms of the flat coordinate system.

\medskip

\noindent To proceed further in the study of these natural submanifolds
one must study the curvature of the induced metric and its relation
with the various other structures.

\medskip

A powerful in the study of Frobenius manifolds is the extended
connection
$\widetilde\nabla$ on $M\times \mathbf{P}^1$ defined by
\begin{eqnarray*}
{\widetilde\nabla}_U & = & \nabla_U + z U \circ \,, \\
{\widetilde\nabla}_{z\frac{d~}{dz}} & = & z\frac{d~}{dz} +
z E\circ \, - \, \nu\,,
\end{eqnarray*}
where $\nu(U) = \frac{D}{2}\, U - \nabla_U E\,,$ (so
$<X,\nu(Y)>+<\nu(X),Y>=0\,.$) The vanishing of the curvature of this extended
connection is then equivalent to the above definition of a
Frobenius manifold.

Since the pull-back of any totally symmetric $(r,0)$ tensor from
$M$ to a submanifold $N$ remains totally symmetric, any natural
submanifold of a Frobenius manifold retains the quasi-potentiality
condition, even if it is not flat (note, however, that there is no reason for an arbitrary
$\F$ manifold to be quasi-potential). Similarly, one may restrict
the extended connection on $M\times\mathbf{P}^1$ to
$N\times\mathbf{P}^1\,.$ This new connection will not be flat, but
it still has special curvature properties.

\begin{prop} The curvature of the restriction of the extended connection of a Frobenius manifold
to a natural submanifold is independent of $z\,.$
\end{prop}

\noindent{\bf Proof}  Consider
\begin{eqnarray*}
{\widetilde R}(U,V)W & = & \left(
{\widetilde\nabla}_U
{\widetilde\nabla}_V-{\widetilde\nabla}_V{\widetilde\nabla}_U -
{\widetilde\nabla}_{[U,V]} \right)W \,, \\
& = & \{ R(U,V) W\}\\
& & + z \{ U\circ \nabla_V W - V \circ \nabla_U W -
[U,V]\circ W
+ \nabla_U (V\circ W) - \nabla_V (U\circ W) \}\\
& & + z^2 \{ U \circ (V
\circ W) - V \circ (U \circ W) \}\,.
\end{eqnarray*}
The last term vanishes by the commutativity and associativity of
the $\circ$ product. The middle term vanishes by a result of
Hertling (using the definition (\ref{Fmandef}) of an $F$ manifold.
the first term is just the curvature of the manifold $N\,.$ The
only
other curvature to calculate is
\begin{eqnarray*}
{\widetilde R}\left( z \frac{d~}{dz},U\right)V & = &
{\phantom{-}}\left(
{\widetilde\nabla}_{z \frac{d~}{dz}}
{\widetilde\nabla}_U -
{\widetilde\nabla}_U
{\widetilde\nabla}_{z \frac{d~}{dz}}
\right) V \,, \\
& = & {\phantom{-}}\nabla_U \,\mu(V) - \mu (\nabla_U V)\,,\\
& = & -\nabla_U \nabla_V E + \nabla_{\nabla_U V} E \,,\\
& = & - \Delta(U,V)\,,
\end{eqnarray*}
the other terms vanishing for similar reasons as above.
Thus the only non-zero terms in the curvature of the extended
connection is the curvature $R(U,V)W$ of the manifold $N$ and
$\Delta(U,V)\,.$ However since
\begin{eqnarray*}
\Delta(U,V)-\Delta(V,U) & = & \nabla_U\nabla_V E - \nabla_V\nabla_U
E - \nabla_{\nabla_U V - \nabla_V U} E \,, \\
& = & \nabla_U\nabla_V E - \nabla_V\nabla_U E - \nabla_{[U,V]}\,,\\
& = & R(U,V)E
\end{eqnarray*}
the independent non-zero terms are curvature and the symmetric part
$\Delta^s(U,V)$ of $\Delta(U,V)\,.$

\bigskip

\medskip

\noindent The following result relates $\Delta$ to the curvature, $\Delta(U,V)=R(U,E)V\,,$
This properties may be proved using submanifold theory. The following theorem is
standard, and is included here only to fix notation.

\medskip

\begin{theorem} Let $M$ be a manifold with Levi-Civita connection
$\overline{\nabla}$ and let $N$ be an arbitrary submanifold. Then
for all $W\,,X\,,Y\,,Z\in TN$ and normal vectors $\xi\,,\eta\in
TN^\perp\,:$
\begin{itemize}
\item{Gauss formula:}
\[
{\overline\nabla}_X Y = \underbrace{ \nabla_X  Y}_{TN} +
\underbrace{ \alpha(X,Y) }_{TN^\perp}\,;
\]
\item{Weingarten formula:}
\[
{\overline\nabla}_X \xi = \underbrace{- A_\xi X}_{TN} +
\underbrace{ \nabla^\perp _X \xi}_{TN^\perp}\,;
\]
\item{Gauss equation:}
\[
<{\overline R} (X,Y)Z,W> = <R(X,Y)Z,W> -
<\alpha(Y,Z),\alpha(X,W)>+<\alpha(X,Z),\alpha(Y,W)>\,;
\]
\item{Codazzi equation}
\[
({\overline R}(X,Y),Z)^\perp =
(\nabla^\perp_X \alpha)(Y,Z)-
(\nabla^\perp_Y \alpha)(X,Z)\,;
\]
\item{Ricci equation:}
\[
<{\overline R}(X,Y) \xi, \eta> = < R^\perp(X,Y) \xi, \eta> -
<[A_\xi, A_\eta]X,Y>\,.
\]
\end{itemize}
Here $\alpha$ is the second fundamental form and $A$ is the
shape operator, which are related by
\[
<\alpha(X,Y),\xi>=<A_\xi X,Y> \quad\quad \forall X,Y \in TN\,, \xi
\in TM^\perp\,.
\]
\end{theorem}

\bigskip

\begin{prop} Let $N$ be a natural submanifold of a Frobenius
manifold $M\,.$ Then $N$ (with the naturally induced structures) is
an $\F$ manifold with, for all $U\,,V \in TN\,:$
\begin{eqnarray*}
\Delta(U,V) & = & R(U,E)V \,,\\
Lie_E \alpha & = & 0 \,,\\
< \nabla_U e^\top , V> & = & < \alpha(U,V) , e^\perp >\,,
\end{eqnarray*}
\end{prop}

\medskip

\noindent{\bf Proof}  The proof is a simple exercise in
submanifold theory. Recall that for a natural submanifold
$E^\perp = 0\,$ so $E\in TN$ and that $\alpha(U,V) \in TN^\perp$
for all $U\,,V \in TN\,.$
\begin{equation}
\begin{array}{rcl}
{\overline\Delta}(U,V) & = & {\overline \nabla}_U {\overline
\nabla}_V E - {\overline\nabla}_{{\overline\nabla}_U V} E \,, \\
& = & \Delta(U,V) - A_{\alpha(V,E)} U + A_{\alpha(U,V)} E\\
&&+\alpha(U, \nabla_V E) - \alpha(\nabla_U V,E) - [\alpha(U,V), E] \\
&&+ \nabla^\perp_U
\alpha(V,E)-\nabla_E^\perp \alpha(U,V)\,,
\end{array}
\label{decomp}
\end{equation}
where the torsion free condition has been used to to calculate the
term ${\overline\nabla}_{\alpha(U,V)} E\,.$ Taking the tangential
component of (\ref{decomp}) yields (since ${\overline\Delta}=0\,$)
\[
\Delta(U,V) =A_{\alpha(V,E)} U - A_{\alpha(U,V)} E+
[\alpha(U,V),E]^\top
\]
and taking the inner product of this with $W\in TN$ gives
\[
<\Delta(U,V) W> = < \alpha(V,E),\alpha(U,W)>-
<\alpha(U,V),\alpha(E,W)> + <[\alpha(U,V),E]^\top,W>\,.
\]
Using (\ref{Fmandef2}) gives $<[\alpha(U,V),E],W>=0$ so
\begin{equation}
[E,\alpha(U,V)]^\top=0\,.
\label{perprel}
\end{equation}
Hence, by the Gauss equation
\begin{eqnarray*}
<\Delta(U,V),W> & = & < R(U,E)V,W> \,,\\
& = & R(U,E,V,W)
\end{eqnarray*}
or $\Delta(U,V)=R(U,E)V\,.$

\medskip

\noindent Taking the perpendicular component of (\ref{decomp}) yields
\[
0 =
\nabla^\perp_U
\alpha(V,E)-\nabla_E^\perp \alpha(U,V)+\alpha(U, \nabla_V E) -
\alpha(\nabla_U V, E) - [\alpha(U,V),E]^\perp\,.
\]
Using the Codazzi equation and the torsion free property of the
induced connection gives
\begin{eqnarray*}
\left(Lie_E \alpha\right) (U,V)  & = &
\left[ Lie_E \alpha(U,V)\right]^\top\,,\\
& = & 0
\end{eqnarray*}
by (\ref{perprel}), so $Lie_E \alpha = 0\,.$

\medskip

\noindent To obtain the last part of the proposition, recall that
the identity field $e$ satisfies the relation ${\overline\nabla}
e=0$ and decomposes as $e=e^\top+ e^\perp$ on $TN\,$ (this decomposition implies that
$e^\top \circ e^\top=e^\top\,, e^\perp \circ e^\perp = e^\perp$ and $e^\top \circ e^\perp=0$
which gives a decomposition of a cohomogeneity one submanifolds). Thus
\begin{eqnarray*}
0 & = &{\overline\nabla}_U e \,, \\
& = & {\overline\nabla}_U e^\top + {\overline\nabla}_U e^\perp \,,\\
& = & \left(\nabla_U e^\top - A_{e^\perp} U\right) + \left( \nabla_U^\perp e^\perp
+ \alpha(U,e^\top) \right)\,.
\end{eqnarray*}
Decomposing this into tangential and perpendicular components
gives
\begin{eqnarray*}
<\nabla_U e^\top,V> & = & < A_{e^\top} U,V> \,, \\
& = & <\alpha(U,V), e^\perp> \,, \\
\nabla^\perp_U e^\perp & = & - \alpha(U, e^\top)\,.
\end{eqnarray*}
[Note that $<\nabla_U e^\top,V>-<\nabla_V e^\top,U>=0$ so
the 1-form $\varepsilon_N(\cdot)=<,e^\top,\cdot>$ is closed, as it
must, since $d (\imath^\star \varepsilon)=0\,.$]

\bigskip

An immediate corollary of this proposition is the following:

\begin{corollary}
Any flat caustic of a semi-simple Frobenius manifold is itself a Frobenius
manifolds, i.e. it is a Frobenius submanifold. All two dimensional
caustics are Frobenius submanifolds.
\end{corollary}

\noindent{\bf Proof} The only thing to note is that for natural
submanifold of a semi-simple Frobenius manifold
\[
\{e^\perp=0\}   \Longleftrightarrow \{ N~{\rm is~a~caustic~}\}\,.
\]
Hence by the above proposition $\nabla_U e=0\,.$ If the caustic is
flat then all obstruction vanish. Note that on a general
(non-flat) caustic $\alpha(e,U)=0$ so from the Gauss equation
\begin{equation}
R(W,X,Y,Z)=0\quad\quad {\rm if~any~of~the~vector~fields~}
W,X,Y,Z=e\,.
\label{curvcaustic}
\end{equation}
Hence all two-dimensional caustics are Frobenius submanifolds.

\bigskip

\subsection{Semi-simple $\F$-manifolds}

In this subsection the curvature properties of semi-simple
$\F$-manifolds will be studied. Again the approach will stress
those properties intrinsic to an $\F$-manifold as defined above,
and those which a natural submanifold of a Frobenius manifold
possesses. From the semi-simplicity and the compatibility of the
multiplication with the metric:
\begin{eqnarray*}
\eta_{ij} & = & < \partial_i, \partial_j > \,, \\
& = & < e , \partial_i \circ \partial_j > \,, \\
& = & \delta_{ij} < e \,, \partial_i>
\end{eqnarray*}
and hence the metric is diagonal. Curvature calculations for
diagonal metrics are standard. For completeness and to fix
notation:

\begin{prop}\label{curvatureprop} Let
$\eta_{ii}=H_i^2\,,\beta_{ij}=\frac{\partial_i H_j}{H_i}\,.$ Then:
\begin{eqnarray*}
\Gamma_{jk}^i  & = & 0\,, {\rm~where~}i,j,k
{\rm~are~distinct}\,;\\
\Gamma_{ik}^i & = &\frac{H_k}{H_i} \beta_{ki}\,;\\
\Gamma^i_{jj} & = &- \frac{H_j}{H_i} \beta_{ij}\,,{\rm~where~} i\neq j\,.
\end{eqnarray*}
Similar formulae hold for $\Gamma^{ij}_k\,,$ where
$\Gamma^{ij}_k=-g^{is} \Gamma^j_{sk}\,.$
Moreover: $R^{ij}_{kl}=0$ if $i,j,k,l$ are distinct,
$R^{ii}_{kl} =
R^{ij}_{kk}=0\,,R^{ij}_{il}=-R^{ij}_{li}=R^{ji}_{li}=-R^{ji}_{il}\,,$
and, for $i\neq j, i\neq l\,,$
\begin{eqnarray*}
R^{ij}_{il} & = &\frac{1}{H_i H_j} \left\{ \partial_l
\beta_{ji}-\beta_{jl} \beta_{li}\right\}\,,\\
R^{ij}_{ij} & = & \frac{1}{H_i H_j} \left\{
\partial_i \beta_{ij} + \partial_j \beta_{ji} +
\sum_{p\neq i,j} \beta_{pj} \beta_{pi} \right\}\,.
\end{eqnarray*}
Here $R^{ij}_{kl} = \eta^{is} R^j_{skl}$ and
$R(\partial_i,\partial_j)\partial_k= R^r_{kij} \partial_r\,.$
\end{prop}

\medskip

\noindent In addition to the metric and multiplication on $\F$ one
has unity and Euler vector fields:
\begin{eqnarray*}
e & = & \sum_{i=1}^n \partial_i \,, \\
E & = & \sum_{i=1}^n u^i \partial_i\,
\end{eqnarray*}
and the homogeneity condition (\ref{Fmandef2})
becomes
$E(\eta_{ii}) = (D-2) \eta_{ii}\,,$ or
$E(\beta_{ij})=  - \beta_{ij}\,.$
No more can be said about the curvature properties of a general $\F$ manifold.

Consider now a semi-simple Frobenius manifold. Condition (b) in
its definition implies:
\[
(b^\prime) \quad\quad \Rightarrow \quad\quad
\{{\rm metric~is~Egoroff,~i.e.~}\beta_{ij}=\beta_{ji}\}\,;
\]
and hence the induced metric on the natural submanifold is
Egoroff, as condition $(b)$ hold for any submanifold (and
it is clear that the induced metric on a natural
submanifold remains diagonal). Thus natural submanifolds
of semi-simple Frobenius manifolds are Egoroff (one may also prove
this directly from the definition of a natural submanifold).

The following proposition may be derived by direct computation in
diagonal coordinates, so no proof will be given.
\begin{prop} Consider a semi-simple Egoroff $\F$ manifold in
canonical coordinates. Then:
\begin{eqnarray*}
\nabla_i \nabla_j E^k = 0 & \Longleftrightarrow &
\partial_k \beta_{ij} - \beta_{ik} \beta_{kj} \,, \quad i,j,k,{\rm
~distinct}\,,\\
\nabla_i \nabla_j E^i = 0 & \Longleftrightarrow &
e(\beta_{ij})=0\,, \quad i\neq j\,.
\end{eqnarray*}
Hence
\[
\Delta=0 \Longleftrightarrow R=0
\]
and so there is only one obstruction to the extended connection on
$N\times\mathbf{P}^1$ being zero, namely $\Delta(U,V)\,.$
\end{prop}

\noindent For an Egoroff metric, the curvature components may be combined
to give
\[
R^{ij}_{ij}+\sum_{p\neq i,j} R^{ij}_{ip} = \frac{1}{H_i H_j}
e(\beta_{ij}) \,, \quad\quad i\neq j\,.
\]
This formula enables one to consider the curvature of pure
caustics and pure disciminants. On a pure caustic
$e(H_i)=0 \Rightarrow e(\beta_{ij})=0$ and this right hand side vanishes (in
accordance with (\ref{curvcaustic})). On a pure discriminant
$\partial_k \beta_{ij}-\beta_{ik}\beta_{kj}=0\,,i,j,k$ distinct
(this will be shown in the next section), so the only non-zero
curvatures are
$R^{ij}_{ij}=\frac{1}{H_i H_j} e(\beta_{ij})\,.$
Pure discriminants also have the property of having flat normal
bundles, $R^\perp=0\,,$ the conditions $\{u^i=0\,,i\in\D\}$ being
a holonomic nets of lines of curvature \cite{F}.

\medskip

\subsection{The induced intersection form and pencils of
compatible metrics}

An important feature of a Frobenius manifold $M$ is the existence of a
second flat metric, the intersection form, defined by
\[
{}^{(2)} g^{ij} = E(dt^i \circ dt^j)
\]
(in what follows the original metric $g$ will be denoted ${}^{(1)}
g$). Here $\{t^i\}$ are the flat coordinates on $M\,.$ One
important feature of this metric is that
\[
Lie_e\,{}^{(2)}\!g^{ij} = {}^{(1)}\!g^{ij}
\]
and it follows from this that the pencil of metric defined by
${}^{(\Lambda)}\! g^{ij} ={}^{(2)}g^{ij}+\Lambda {}^{(1)}g^{ij}$
is flat for all values of $\Lambda\,.$

A second metric, and hence a pencil of (inverse)-metrics, may also
be defined on a natural submanifold $N\subset M$ since, by
definition, $TN\circ TN\subset TN$ and $E_x \in TN\,\forall x\in
TN\,.$ This is not immediately obvious for a discriminant
submanifold since an equivalent definition of canonical
coordinates for a Frobenius manifold is as solutions of the
polynomial equation
\begin{equation}
\det[ {}^{(2)}\! g^{ij} - u \,{}^{(1)}\! g^{ij} ] = 0\,,
\label{canonicaldef}
\end{equation}
and hence on a submanifold with $u^i=0\,,\det[ {}^{(2)}\!
g^{ij}]=0$ so the metric is non-invertible. However the problem
lies in the orthogonal component to $N\,.$ One may use the
orthogonal decomposition $TM\cong TN \oplus TN^\perp\,,$ and
consider
\[
g^{ij} \frac{\partial~}{\partial t^i} \otimes_s
\frac{\partial~}{\partial t^j} \in
TM \otimes_s TM  \cong (TN \otimes_s TN) \oplus (TN \otimes_s
TN^\perp) \oplus (TN^\perp \otimes_s TN^\perp)\,.
\]
The cross term vanishes and hence on obtains a symmetric bilinear
from on $TN \otimes_s TN\,.$  In canonical coordinates (and hence
for a semi-simple Frobenius manifold)
\[
{}^{(2)} g_{ij} = \sum_{i=1}^m \frac{\eta_{ii}}{u^i} (du^i)^2\,,
\]
and similar looking formulae hold for the induced metrics on a
natural submanifold. Although both induced metric will no longer,
in general, be flat, certain special curvature properties remain.

\begin{lemma}\label{curvaturelemma} Consider two diagonal metrics $(i=1,2)$
\[
{}^{(i)}\! g = \sum_{r=1}^n {}^{(i)}\! g_{rr} (du^r)^2
\]
with rotation coefficients ${}^{(r)}\! \beta_{ij}$ and
${}^{(i)}\!H_r=\sqrt{{}^{(i)}\! g_{rr}}\,.$ Let
\[
{}^{(2)}\! g_{rr}=\frac{{}^{(1)}\! g_{rr}}{u^r}
\]
and
\begin{equation}
{}^{(\Lambda)}\! {g}_{rr} ={}^{(2)}\! g_{rr} +\Lambda\, \, {}^{(1)}\! g_{rr}\,.
\label{pencil}
\end{equation}
Then
\begin{eqnarray*}
{}^{(\Lambda)}\! \Gamma^{ij}_k & = & {}^{(2)}\! \Gamma^{ij}_k+
\Lambda\,\, {}^{(1)}\! \Gamma^{ij}_k
\,, \\
{}^{(\Lambda)}\! R^{ij}_{rs} & = &
{}^{(2)}\! R^{ij}_{rs}+ \Lambda \,\,{}^{(1)}\! R^{ij}_{rs}
\end{eqnarray*}
where ${}^{(\Lambda)}\! \Gamma^{ij}_k$ and ${}^{(\Lambda)}\! R^{ij}_{rs}$
are the appropriate Christoffel and curvatures of the pencil of
metric (\ref{pencil})
\end{lemma}
Thus any semi-simple $\F$ manifold carries such a pencil, and in
particular, so does any natural submanifold of a Frobenius
manifold. In the terminology of \cite{Mo}, on such (sub)-manifolds one
has a pencil of compatible metrics. It is this result that will be
behind the study of induced bi-Hamiltonian structures on natural
submanifolds that will be given in the next section. Pure
disciminant submanifold have further special properties:

\begin{theorem}[\cite{D}] Let $M$ be a Frobenius manifold and let
$N$ be a pure discriminant submanifold. Then the metric on $N$
induced from ${}^{(2)}g$ is flat.
\end{theorem}

\noindent Thus on a pure discriminant one has a distinguished
coordinate system being the flat coordinates for the second
induced metric. Examples will be given in the next section (see
also Main Example B). Since
\[
{}^{(2)}\! \beta_{ij} = \sqrt{\frac{u_i}{u_j}} \left( {}^{(1)}\!
\beta_{ij}\right)
\]
(and note that the second metric will not be Egoroff) it follows
from
\[
\partial_k {}^{(2)} \!\beta_{ij}-{}^{(2)}\! \beta_{ik}{}^{(2)}
\!\beta_{kj}= \sqrt{\frac{u_i}{u_j}}\left(
\partial_k {}^{(1)} \!\beta_{ij}-{}^{(1)} \!\beta_{ik}{}^{(1)} \!\beta_{kj}
\right)
\]
that ${}^{(1)} \!R^{ij}_{ik}=0\,,$ as stated
at the end of the last section.

\subsection{Tangent vectors to a natural submanifold}

The intersection form of a Frobenius manifold $M$ enables one to
construct natural vector fields tangent to a natural submanifold.
In a flat coordinate system $\{t^i\}$ where $e=\partial_1\,,$
${}^{(2)} g^{mi} = E^i\,,$ so the components of the last
row/column of the intersection form are the components of the
Euler vector field which, by definition, is tangent to a natural
submanifold. Consider the vector fields defined on $N$ by
\[
V^{(\alpha)} = \left.\left( \frac{\partial
t^i}{\partial\tau^\alpha} {}^{(1)} \!g_{ij} {}^{(2)}\! g^{jk} \right)
\right|_N \frac{\partial~}{\partial t^k}
\]
which a priori lie in $TM\,.$ Using (\ref{orthogdecomp}) and the
fact that
\[
\left.E^k\right|_N = E_N^\alpha \frac{\partial
t^k}{\partial\tau^\alpha}
\]
it follows that the component of $V^{(\alpha)}$ in $TN^\perp$ is
\[
\left(\Xi_{\alpha\beta}^{\tilde\alpha}
E_N^\beta\right)\frac{\partial~}{\partial \nu^{\tilde\alpha}}
\]
and hence is zero.

\medskip

Using techniques identical to the above, one may show that if
\begin{equation}
X\circ Y \in TN \quad\quad \forall X\in TN\,, \quad Y\in TM
\label{discriminantcondition}
\end{equation}
then the $m$ vector fields
\[
V^{(i)} = \,{}^{(2)}\!g^{ij} \frac{\partial~}{\partial t^j}
\]
are all tangential to the submanifold. In the semi-simple case
one may show, using their explicit parametrization, that the
only natural submanifolds which satisfy the
above condition are pure discriminants. For Frobenius manifolds
based on Coxeter groups this result is already known \cite{A,Sh}\,.
For codimension one discriminants a simple proof that (\ref{discriminantcondition})
holds may be given using the decomposition of the unity vector field
$e=e^\top+e^\perp\,.$ Since, if $N$ is a natural submanifold,
$X\circ e^\perp=0\,, \forall X\in TN$ then
because any vector in $TN^\perp$ must be
a multiple of $e^\perp$ the result follows.

\subsection{Examples of Frobenius submanifolds}

Frobenius submanifolds certainly exist - any flat caustic, if they
exist, of a semi-simple Frobenius manifold will inherit the
structure of a Frobenius manifold. In theory this gives a way to
find such submanifolds, though in practice it would be
computationally difficult. A more practical way is to look for
submanifolds which are hyperplanes (in the flat coordinates
$\{t^i\}$, and coordinate hyperplanes in particular.

\medskip

\begin{example} {\rm Let $I\subset\{1\,,2\,,\ldots\,,m\}\,$ and
suppose that $\N$ is given by the conditions $t^i=0$ for $i\notin
I\,.$ Then the obstruction reduces to the algebraic condition
\[
\left.c_{ij}^{~~~k}\right|_N = 0 \,,\quad\quad i\,,j \in I\,,
k\notin I\,.
\]
This condition was derived in \cite{Z} in the context of Frobenius
manifolds constructed from Coxeter groups. Here it is a
specialization of the more general condition (\ref{obstruction}).}
\end{example}

\medskip

\begin{example} {\rm \cite{M}, section III.8.7.1.}\end{example}

\medskip

\noindent It is not clear, but follows from the results above,
that these submanifolds are caustics.

\medskip

Large numbers of examples may be found using these results. For
example, for Frobenius manifolds constructed from Coxeter groups
one finds that the submanifold associated with another Coxeter
group obtained by \lq folding\rq~the original Coxter diagram. For
example the $H_3$-Frobenius contains the Frobenius manifold
$I_2(10)$ as a submanifold, which corresponds to the folding

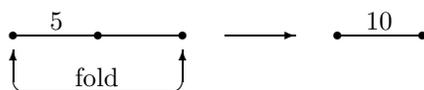
\setlength{\unitlength}{0.75cm}
\begin{picture}(11,2)
\put(4,1){\circle*{0.15}} \put(4,1){\line(1,0){1.5}}
\put(5.5,1){\circle*{0.15}} \put(5.5,1){\line(1,0){1.5}}
\put(7,1){\circle*{0.15}} \put(7.75,1){\vector(1,0){1.25}}
\put(9.75,1){\circle*{0.15}} \put(9.75,1){\line(1,0){1.5}}
\put(11.25,1){\circle*{0.15}} \put(4.65,1.1){5}
\put(10.25,1.1){10}

\put(5.5,0.25){\oval(3,0.5)[b]}

\put(4,0.25){\vector(0,1){0.5}} \put(7,0.25){\vector(0,1){0.5}}

\put(5.1,0.10){fold}

\end{picture}

\vskip 2mm

\noindent The possible foldings, and hence submanifolds, of Coxter
groups, are given in Table 1.

\bigskip

\begin{table}
\begin{center}
\begin{tabular}{c|c}
\hline
Coxeter Group & Coxeter subgroup \\ \hline
$A_{2n+1}$ & $B_n$ \\
$D_{n+1}$ & $B_n$ \\
$D_5$ & $H_3$\\
$E_6$ & $F_4$\\
$E_8$ & $H_4$\\
$W$ (arbitrary) & $ I_2({\rm~Coxeter~number~of~}W)$\\
\hline
\end{tabular}
\end{center}
\vskip 5mm
\caption{Frobenius submanifolds associated with Coxeter groups}
\end{table}

\vskip 5mm

\noindent  These results may also be generalized to Frobenius
manifolds constructed from extended affine Coxeter groups \cite{DZ1}(and
probably more generally too, since notions of foldings exist more
generally in singularity theory). The analogous results are given in Table 2.

\medskip

\begin{table}
\begin{center}
\begin{tabular}{c|c}
\hline
Extended affine Coxeter Group & Extended affine Coxeter subgroup \\ \hline
$A_{l=2n-1}^{k=n}$ & $C_n$ \\
$D_{n+1}$ & $B_n$ \\
$E_6$ & $F_4$\\
$W$ (arbitrary) & $A^{k=1}_{l=1}$\\
\hline
\end{tabular}
\end{center}
\vskip 5mm
\caption{Extended affine Frobenius submanifolds}
\end{table}

\vskip 5mm

However, such submanifolds are just hyperplanes. More interesting examples
may be constructed.

\medskip

\begin{example} \label{A3example} {\rm Frobenius submanifolds of $A_3$
(this is a special case of Main Example A). For the
$A_3$-singularity one takes the polynomial
\[
p(z)=z^4 + a_1 z^2 + a_2 z + a_3
\]
and constructs the metric via the formula (\ref{Anmetric}).
This metric is flat, though not in flat coordinates, which are given by
\begin{eqnarray*}
a_1 & = & t_3 \,, \\
a_2 & = & t_2 \,, \\
a_3 & = & t_1 + \frac{1}{8} t_3^2\,.
\end{eqnarray*}
In these coordinates the metric takes the standard antidiagonal form.
Having fixed the flat coordinates, so
\begin{equation}
p(z)=x^4 + t_1 + \frac{1}{8} t_3^2 + t_2 x + t_3 x^2
\label{deff}
\end{equation}
the algebra is defined by
\[
c_{ijk} = - \res_{x=\infty}
\frac{
\partial_{t_i} f \,
\partial_{t_j} f \,
\partial_{t_k} f }{\partial_z f}\, dz
\,.
\]
From this the prepotential may be constructed. The canonical
coordinates are now defined as roots of the cubic
(\ref{canonicaldef}). The discriminant and caustics may
easily be calculated from this cubic.

\bigskip

\noindent{\underline{Discriminant}} (The induced structure on the
disciminant does not define a Frobenius manifold. The results are
given here for use in the next section). If $u^i=0$ then from
(\ref{canonicaldef}) $\det({}^{(2)}\!g^{ij})=0$, or
\begin{equation}
-t_1^3 + \frac{27}{256} t_2^4 - \frac{9}{16} t_1 t_2^2 t_3
+ \frac{1}{8} t_1^2 t_3^2 - \frac{7}{128} t_2^2 t_3^3 +
\frac{1}{64} t_1 t_3^4 - \frac{1}{512} t_3^6=0\,.
\label{discrim}
\end{equation}
This is precisely the condition for the polynomial (\ref{deff}) to
have a repeated root, i.e. it defines a discriminant hypersurface.
Using the fact that such surfaces are ruled one may easily obtain
a parametrization of the surface
\begin{eqnarray*}
t_1 & = & + 2^{-9}( u^4 - 6 u^2 v^2 - v^4)\,,\\
t_2 & = & + 2^{-4} u v^2 \,, \\
t_3 & = & - 2^{-3}(u^2 + v^2)\,,
\end{eqnarray*}
where $u$ and $v$ are the flat coordinates for the induced intersection form.
Such an explicit parametrization will be used in the next section to construct
induced
bi-Hamiltonian structures on this discriminant.

\bigskip

\noindent{\underline{Caustics}} If $u^i=u^j$ for some $ i \neq j$
then the
polynomial (\ref{canonicaldef}) must have a double root, and the
condition for this is either (a) $t_2=0$ or
(b) $ 27 t_2^2 + 8 t_3^3=0\,.$ These surfaces correspond to the
cylinders over the Maxwell strata  and caustic of the polynomial
(\ref{deff}). The induced structures on these surfaces define
Frobenius manifolds (corresponding to the Coxeter group $I_2(4)$):
the structure on the (a) being studied by Zuber \cite{Z} and on (b) by
the author \cite{S2}. Note how the two parts in the definition
of a natural submanifold are used; the induced multiplication on the
flat surface $t_2^2 + k t_3^3=0$ is associative for all values of $k$, but
the Euler field is tangential only for the two special values of $k$ given above.}

\end{example}

\bigskip

Further examples may be obtained by tensoring Frobenius manifolds
together and restricting structures to various hyperplanes \cite{M,S2}.

\bigskip

\section{Induced bi-Hamiltonian structure}

An equation of hydrodynamic type is, by definition, of the form
\begin{equation}
U^i_T=V^i_j(U^k) U^j_X\,.
\label{hydro}
\end{equation}
It was observed by Riemann that
such system transform covariantly with respect to arbitrary
changes of dependent variables ${\tilde U}^i = {\tilde U}^i(U)\,.$
It is not surprising therefore that geometrical ideas should be used in
the study of such equations. Associated with any semi-simple Frobenius manifold
is a bi-Hamiltonian hierarchy of such hydrodynamic equations. The aim of
this section is to show how one may constrain such systems onto a
natural submanifold while retaining its bi-Hamiltonian structure.
In this case the system (\ref{hydro}) may be diagonalised
\[
u^i_T = \lambda^i(u) u^i_X\,,\quad\quad i=1\,,\ldots\,,m\,.
\]
(Note that a diagonalised system is defined in terms of Riemann
invariants
\[
R^i_T = \lambda^i(R) R^i_X
\]
and these are only defined up to $R^i\mapsto{\tilde R}^i(R^i)$
transformation. The canonical coordinates are specific examples
of Riemann invariants, and do no have such a freedom.) The $\lambda^i$ are known as
the characteristic speeds of the system. It is
immediate from this diagonal form that if the mild finiteness
condition
\begin{equation}
\left.\lambda^j\right|_{\lambda^i=0} < \infty
\label{conditionD}
\end{equation}
holds then one may restrict the system to the discriminant
$u^i=0\,.$ To reduce the system to a caustic $u^i-u^j=0$ requires the much
stronger condition
\begin{equation}
\left. (\lambda^i-\lambda^j)\right|_{u^i-u^j=0} = 0
\label{conditionC}
\end{equation}
together with a mild finiteness condition for the remaining
${\lambda^i}\,.$ Such constraints do hold for the systems
associated with semi-simple Frobenius manifolds and this will be proved below.

\medskip

\noindent Before this a more general discussion will be given which will
associate $\F$ manifolds with certain diagonal sets of equations.

\subsection{Semi-Hamiltonian systems and curved $\F$-manifolds}

\begin{Def} A diagonal system of hydrodynamic type
\[
u^i_T = \lambda^i(u) u^i_X\,,\quad\quad i=1\,,\ldots\,,m\,,
\]
is semi-Hamiltonian if there exists a diagonal metric
\[
g=\sum_{i=1}^m g_{ii}(u) (du^i)^2
\]
satisfying the equations
\[
\partial_j \log\sqrt{g_{ii}} = \frac{\partial_j \lambda^i}{\lambda^j-\lambda^i}
\]
for $i\neq j\,.$
\end{Def}
On cross differentiating on obtains the identities
\[
\partial_k \frac{\partial_j \lambda^i}{\lambda^j-\lambda^i}=
\partial_j \frac{\partial_k \lambda^i}{\lambda^k-\lambda^i}\,,
\]
for distinct $i,j,k\,.$ All semi-Hamiltonian systems are
integrable, via the generalized hodograph transform. Such
systems possess an infinite number of commuting flows
\[
U^i_{T^\prime} = w^i(U) U^i_X
\]
where the $w^i$ are solutions of the linear system
\[
\frac{\partial w^i}{w^j-w^i} = \frac{\partial_j \lambda^i}{\lambda^j-\lambda^i}\,,
\quad\quad i\neq j,,
\]
which exists since its integrability follows from the
definition of semi-Hamiltonian. If the metric is homogeneous
with respect to the vector field $E$, so
\[
E \eta_{ii} = (D-2) \eta_{ii}
\]
them one obtains a semi-simple $\F$ manifold. It also follows
from the
definition of semi-Hamiltonian that the only non-zero curvature
components are $R^{ij}_{ij}\,,$ for $i\neq j\,.$

\bigskip

\noindent{\bf Example} Consider the system
\begin{equation}
u^i_T = \left( \sum_{r=1}^m u^r + 2 u^i\right) u^i_X\,, \quad\quad i=1\,,\ldots\,, m\,,
\label{multidkdv}
\end{equation}
which corresponds to the dispersionless limit of the coupled KdV hierarchy \cite{FP,AF}.
Such a system is semi-Hamiltonian with metric
\[
g=\sum_{i=1}^m \left\{ \frac{\prod_{r\neq i} (u^r-u^i)}{\phi_i(u^i)}\right\}
(du^i)^2\,, \quad\quad \phi_i {\rm~arbitrary\,.}
\]
The metric ${}^{(1)}\!g$ is defined as the above metric with $\phi_i=1\,,$ and is homogeneous
with respect to the Euler vector field, and so one obtained  semi-simple
$\F$ manifold. This metric is not Egoroff, so one does not obtain a Frobenius
manifold, nor can it be a natural submanifold of a semi-simple Frobenius manifold.
This metric is, however flat; in flat coordinates
\[
{}^{(1)}\!g^{ij} =
\begin{pmatrix}
  0 & 0 & 1 \\
  0 & 1 & 0 \\
  1 & 0 & 0
\end{pmatrix}
\,,\quad\quad\quad
{}^{(2)}\!g^{ij} =
\begin{pmatrix}
  0 & -1 & t_1/2 \\
  -1 & 0 & t_2/2 \\
  t_1/2 & t_2/1 & t_3
\end{pmatrix}\,,
\]
(both flat) and the Euler and identity fields are
\begin{eqnarray*}
E & = & t_1 \partial_1+2 t_2 \partial_2 + 3 t_3 \partial_3\,,\\
e & = & 3 \partial_1 + \frac{1}{2} t_1 \partial_2 - \frac{1}{2} t_2 \partial_3\,.
\end{eqnarray*}
The associative multiplication is somewhat complicated when written in the
flat coordinate system. Note also that
\[
\det[ {}^{(1)}\!g ] = \prod_{i\neq j} (u^i-u^j)^2\,,
\]
so the metric is degenerate on bifurcation diagrams.

\medskip

This example may also be used to illustrate how one may reduce hydrodynamic systems
to bifurcation diagrams and discriminants, and also some of the problems that may arise.
Restricting the system (\ref{multidkdv}) to a disciminant is trivial - the form of the
equations is unchanged. On a bifurcation diagram $\{k_1,\ldots,k_r\}$ one obtains
\[
u^i_T = \left( \sum_{r=1}^n k_r u^r + 2 u^i\right) u^i_X\,, \quad\quad i=1\,,\ldots\,, n\,,
\]
which remains semi-Hamiltonian (in fact for all values for the $k_i$) with metric
\[
{}^{(1)} g_{_\mathcal{C}} =
\sum_{i=1}^n \prod_{r\neq i} (u^r-u^i)^{k_i}
(du^i)^2\,.
\]
Thus one obtains a \lq stratified\rq~space where ${}^{(1)}\!g$ is defined everywhere except
on bifurcation diagrams, but with another metric ${}^{(1)}\! g_{_\mathcal{C}}$ on the bifurcation diagram (which in term
is defined everywhere on the bifurcation diagram except at sub-bifurcation diagrams etc.).

\subsection{BiHamiltonian hierarchies associated to Frobenius manifolds}

Given a Frobenius manifold there exits an associated hierarchy of
hydrodynamics type, which, in terms of flat coordinates
$\{t^\alpha\}$ are
given by
\begin{equation}
\frac{\partial t^\gamma}{T^{(\alpha,p)}} = c_{(\alpha,p)\beta}^{~~~\gamma}
\partial_X t^\beta\,, \quad \quad\quad\quad
\left\{
\begin{matrix}
\alpha & = & 1\,, \ldots\,, m\,,\\
p & = & 0 \,, \ldots\,, \infty\,.
\end{matrix}
\right.
\label{frobhierarchy}
\end{equation}
The primary part of the hierarchy is defined by
\[
\frac{\partial t^\gamma}{T^{(\alpha,0)}} = c_{\alpha\beta}^{~~~\gamma}
\partial_X t^\beta\,.
\]
From the flatness of the extended connection it follows that there
exists functions $h_{(\alpha,p)}$ which satisfy the relation
\[
\frac{\partial^2 h_{(\alpha,p)}}{\partial t^\alpha t^\beta}=
c_{\alpha\beta}^{~~~~\gamma}\frac{\partial h_{(\alpha,p-1)}}{\partial t^\gamma}
\]
with the initial condition $h_{(\alpha,0)}=t_\alpha=\eta_{\alpha\beta} t^\beta\,.$
These define Hamiltonians
\[
H_{(\alpha,p)} = \int h_{(\alpha,p+1)} dX
\]
which are conserved with respect to all flows. Here the Hamiltonian
structure is defined,
by the fundamental theorem of Dubrovin and Novikov, by any flat metric $g$.
For functionals $F=\int
f(t,t_X,\ldots) dX\,,G=\int g(t,t_X,\ldots) dX$ the Hamiltonian
is defined by
\[
\{F,G\}=\int \frac{\delta F}{\delta t^i} A^{ij} \frac{\delta
G}{\delta t^j}\, dX
\]
where
\[
A^{ij} = g^{ij}(t) \frac{d~}{dX} - g^{is} \Gamma^j_{sk}(t)
t^k_X\,.
\]
Here $g^{ij}$ is the (inverse) flat metric and $\Gamma_i^{jk}$ the
corresponding Christoffel symbols. The zero-curvature condition
ensures that the bracket satisfies the Jacobi identity.

\bigskip

On a given Frobenius manifold $M$ one has a pencil of flat
metrics ${}^{(\Lambda)} \!g^{ij}={}^{(2)}\!g^{ij} + \Lambda {}^{(1)}\!g^{ij}\,,$ and this
then gives rise to a bi-Hamiltonian structure
\[
\{F,G\}_\Lambda = \{F,G\}_2 + \Lambda \{F,G\}_1\,.
\]
Equivalently, this
then implies
\[
\{ \cdot , H_{(\alpha,p)} \}_1 = {\rm constant~} \{\cdot , H_{(\alpha,p-1)} \}_2\,
\]
so
\begin{equation}
\frac{\partial t^\beta}{\partial T^{(\alpha,p)}}  =
\{t^\beta, H_{(\alpha,p)}\}_1 = {\rm constant~} \{t^\beta, H_{(\alpha,p-1)}\}_2\,,
\label{biham}
\end{equation}
for some constant which depends on various normalisations.

\bigskip

\begin{theorem} Let $M$ be a semi-simple Frobenius manifold. Then
the
restriction of the bi-Hamiltonian hierarchy (\ref{biham}) to
a natural submanifold remains bi-Hamiltonian.
\end{theorem}

\bigskip

\noindent{\bf Proof } Given an arbitrary submanifold $N$  of a flat manifold $M$ one
may constrain, using the Dirac procedure
the corresponding Hamiltonian structure on $M$ to the submanifold
$N\,$ \cite{F}. This results in a Hamiltonian structure of the form
\[
\{F,G\}=\int \frac{\delta F}{\delta \tau^\alpha} A^{\alpha\beta} \frac{\delta
G}{\delta \tau^\beta}\, dX
\]
where
\[
A^{\alpha\beta} = g^{\alpha\beta}(\tau) \frac{d~}{dX} - g^{\alpha\mu} \Gamma^\beta_{\mu\nu}(\tau)
\tau^\nu_X+ \sum_{\tilde\alpha} w^\alpha_{{\tilde\alpha}\mu} \tau^\mu_x
\left( \nabla^\perp\right)^{-1}
w^\beta_{{\tilde\alpha}\nu}\tau^\nu_x\,.
\]
Here $g^{\alpha\beta}$ is the (inverse) induced metric on
$N\,,\Gamma^\alpha_{\mu\nu}$ the corresponding Christoffel
sysmbols, and $w^\alpha_{{\tilde\alpha}\beta}$ the Weingarten
operators of the submanifold. The operator $\nabla^\perp$ is
defined by
\[
\nabla^\perp\phi_\alpha = \frac{d~}{dX} \phi_\alpha + \omega_\alpha^{~\beta} \phi_\beta
\]
where $\omega_\alpha^{~\beta}$ are the normal connection one-forms.

Recall, from lemma \ref{curvaturelemma}, that on a natural submanifold
\begin{eqnarray*}
{}^{(\Lambda)}\! \Gamma^{ij}_k & = & {}^{(2)}\! \Gamma^{ij}_k+
\Lambda\,\, {}^{(1)}\! \Gamma^{ij}_k
\,, \\
{}^{(\Lambda)}\! R^{ij}_{rs} & = &
{}^{(2)}\! R^{ij}_{rs}+ \Lambda \,\,{}^{(1)}\! R^{ij}_{rs}\,.
\end{eqnarray*}
This then implies that if one restricts the flat bi-Hamiltonian structure
on the Frobenius manifold to any natural submanifold one obtains a new
bi-Hamiltonian structure with non-local tails, as above. Care has to be taken with
the structure of the non-local tail for the curved pencil; it is twice as
long as the codimension, each half containing the non-local tail of one of the
individual metrics.

\endproof

\bigskip

It is not obvious from this result that the resulting system is still local.
This may be shown to be case by directly studying the restriction of the
hierarchy (\ref{frobhierarchy}) onto a natural submanifold.
In order to show this the system will first be rewritten in canonical coordinates\footnote{N.B.
The notation used in the rest of this
section differs from that used above. Greek letters denote components of objects
in flat coordinates, and Latin letters denote components of objects in canonical coordinates. Also
$\eta$ will refer to the metric ${}^{(1)}\!g$.}. Extensive
use will be made of the following formulae, all of which are derived in \cite{D}:
\[
c_{\alpha\beta\gamma} = \sum_{i=1}^m
\frac{\psi_{i\alpha}\psi_{i\beta}\psi_{i\gamma}}{\psi_{1i}}\,,
\]
\[
\frac{\partial t^\alpha}{\partial u^i} = \psi_{i1}{\psi_i^{~\alpha}}\,,\qquad\qquad
\frac{\partial u^i}{\partial t^\alpha} = \frac{\psi_{i\alpha}}{\psi_{i1}}\,,
\]
where $\psi_{i1}^2=\eta_{11}\,.$ Indices on $\psi$ are raised and lowered using
$\eta_{\alpha\beta}$, so $\psi_i^{~\alpha} = \psi_{i\beta}\eta^{\beta\alpha}$
and
\[
\sum_{i=1}^m \psi_{i\alpha} \psi_{i\beta} = \eta_{\alpha\beta}\,,\quad\quad
\sum_{i=1}^m \psi_i^{~\alpha} \psi_i^{~\beta}= \eta^{\alpha\beta}\,.
\]

\noindent and crucially the following:
\begin{equation}
(u^j-u^i) \beta_{ij} = \sum_{\alpha} (q_\alpha-\frac{d}{2}) \psi_{i\alpha}
\psi_j^{~\alpha}\,.
\label{useful}
\end{equation}

In canonical coordinates these become
\[
\nabla_i \nabla_j h_{(\alpha,p)} = \delta_{ij} \frac{ \partial h_{(\alpha,p-1)}}{\partial u^j}
\]
and
\begin{eqnarray}
\lambda^i_{(\alpha,p)} & = & \eta^{ii}\nabla_i \nabla_i h_{(\alpha,p)}\,, \label{eqna}\\
& = & \frac{1}{\eta_{ii}} \frac{\partial h_{(\alpha, p-1)}}{\partial u^i} \label{eqnb}
\end{eqnarray}
where $\nabla$ is the Levi-Civita connection of the first metric, and, via the
bi-Hamiltonian property,
\begin{equation}
\lambda^i_{(\alpha,p)} = c_{(\alpha,p)} \, {}^{(2)}\!g^{ii}
{\widetilde\nabla}_i {\widetilde\nabla}_i h_{(\alpha,p)}\,.
\label{eqnc}
\end{equation}
Here $c_{(\alpha,p)}$ are some normalization constants, the precise form of
which will not be required.

In canonical coordinates the characteristic speeds of the primary part of the
hierarchy are given by
\[
\lambda^i_{(\alpha,0)} = \frac{\partial u^i}{\partial t^\alpha}\,.
\]
Thus the conditions (\ref{conditionD}) and (\ref{conditionC}) become, respectively, the
conditions
\begin{equation}
\left.\frac{\partial u^i}{\partial t^\alpha}\right|_{u^j=0} < \infty\,,
\label{frobconditionD}
\end{equation}
and
\begin{equation}
\left.\frac{\partial (u^i-u^j)}{\partial t^\alpha}\right|_{u^i-u^j=0} = 0\,.
\label{frobconditionC}
\end{equation}
If these conditions hold for the primary part of the hierarchy, then they hold for the
entire hierarchy.

\begin{theorem} Suppose that condition (\ref{frobconditionC}) holds. Then
\[
\left.\left( \lambda^i_{(\alpha,p)} - \lambda^j_{(\alpha,p)}\right)\right|_{u^i-u^j=0} = 0
\]
and hence the entire hierarchy may be restricted onto the bifurcation diagram $u^i-u^j=0\,.$
\end{theorem}

\medskip

\noindent{\bf Proof} On expanding equations (\ref{eqna}) and (\ref{eqnc}) one obtains
two different expressions for $\partial_i^2 h\,:$
\begin{eqnarray*}
\lambda^i_{(\alpha,p)} & = & \frac{1}{\eta_{ii}} \left[
\frac{\partial^2~}{\partial {u^i}^2} - \sum_r {}^{(1)}\!\Gamma^r_{ii} \frac{\partial~}{\partial u^r}
\right] h_{(\alpha,p)}\,, \\
\lambda^i_{(\alpha,p+1)} & = & c_{(\alpha,p)} \frac{u^i}{\eta_{ii}} \left[
\frac{\partial^2~}{\partial {u^i}^2} - \sum_r {{}^{(2)}\!\Gamma}^r_{ii} \frac{\partial~}{\partial u^r}
\right] h_{(\alpha,p)}\,.
\end{eqnarray*}
Hence, on using (\ref{eqnb}),
\[
\lambda^i_{(\alpha,p+1)} - \lambda^i_{(\alpha,p)} c_{(\alpha,p)} u^i =
c_{(\alpha,p)} \frac{u^i}{\eta_{ii}}
\left[ \sum_r ({}^{(1)}\!\Gamma^r_{ii} - {{}^{(2)}\!\Gamma}^r_{ii} ) \eta_{rr} \lambda^r_{(\alpha,p+1)}
\right]\,.
\]
Expanding the Christoffel symbols (see Proposition \ref{curvatureprop}).
\[
\left(1-\frac{c_{(\alpha,p)}}{2}\right) \lambda^i_{(\alpha,p+1)} - c_{(\alpha,p)} u^i \lambda^i_{(\alpha,p)}=
c_{(\alpha,p)} \sum_{r\neq i} \frac{1}{H_i H_r} (u^r-u^i) \beta_{ri} \lambda^r_{(\alpha,p+1)}
\]
and on using (\ref{useful})
\begin{eqnarray*}
\left(1-\frac{c_{(\alpha,p)}}{2}\right) \lambda^i_{(\alpha,p+1)} -
c_{(\alpha,p)} u^i \lambda^i_{(\alpha,p)}
&=&
c_{(\alpha,p)} \sum_{\beta, r\neq i}
\left( q_\beta - \frac{d}{2} \right) \frac{\psi_r^{~\beta}}{\psi_{ri}}
\frac{\psi_{i\beta}}{\psi_{i1}}
\lambda^r_{(\alpha,p+1)}\,, \\
&=&
c_{(\alpha,p)} \sum_{\beta, r\neq i}
\left[\left( q_\beta - \frac{d}{2} \right) \frac{\psi_r^{~\beta}}{\psi_{ri}}
\lambda^r_{(\alpha,p+1)}\right] \lambda^i_{(\beta,0)}\,,\\
&=&
c_{(\alpha,p)} \sum_{\beta,\gamma, r\neq i}
\left[\left( q_\beta - \frac{d}{2} \right) \lambda^r_{(\gamma,0)}
\lambda^r_{(\alpha,p+1)}\right] \lambda^i_{(\beta,0)}\,,
\end{eqnarray*}
since
\[
\lambda^i_{(\beta,0)} = \frac{\partial u^i}{\partial t^\beta} = \frac{\psi_{i\beta}}{\psi_{i1}}\,.
\]
This shows two things; firstly that if one can restrict the primary part of the
hierarchy onto a discriminant, then one may restrict the entire hierarchy, and also
that
\[
\left.\left( \lambda^i_{(\alpha,p)} - \lambda^j_{(\alpha,p)}\right)\right|_{u^i-u^j=0} = 0
\]
if the results holds for $\alpha=0\,.$ Hence the result.

\endproof

\noindent As already pointed out, these results are formal, depending on the nondegeneracies
of various quantities when restricted to the submanifold. An example of what can
happen when degeneracies appear will be given in the next section, where the results are
illustrated by means of various examples.

\bigskip

\section{Examples}

The main examples A and B in the introduction are concrete examples of the
general theory developed in this paper. The following examples are
self-explanatory. Further examples may be found in \cite{S2}.

\begin{example} {\rm {\bf The $A_n$ caustics}

The construction of the Frobenius manifold based on the Coxeter group
$A_n$ was given in
Main Example A. It follows that
\[
\frac{\partial u^i}{\partial t^\alpha}=
\left.\frac{\partial p}{\partial t^\alpha}\right|_{z=\alpha^i}\,.
\]
On a caustic $u^i=u^j$ implies that $\alpha^i=\alpha^j$ (which is not
true on a Maxwell strata), so
\[
\left.\frac{\partial u^i}{\partial t^\alpha}\right|_{u^i-u^j=0}=
\left.\frac{\partial p}{\partial t^\alpha}\right|_{z=\alpha^i\,, \,\alpha^i=\alpha^j}=
\left.\frac{\partial u^j}{\partial t^\alpha}\right|_{u^i-u^j=0}\,.
\]
Hence one may restrict the $A_n$ hierarchy onto any caustic. The same
should be true for restriction onto Maxwell strata.}
\end{example}

\medskip

\begin{example} {\rm {\bf  The $A_3$ discriminant}

This example is a continuation of example \ref{A3example}.
Using the parametrization of the swallowtail discriminant given there, the induced metrics
become
\begin{eqnarray*}
\eta_N & = & (-u^4 + 3 u^2 v^2 + v^4)\, du^2 + 2 u v ( u^2 + 4 v^2)
\,du\,dv + v^2 (7 u^2 + v^2) \, dv^2\,,\\
g_N & = & - du^2 - dv^2 \,.
\end{eqnarray*}
Note that $\det \eta_N=v^2(v^2-2 u^2)^3\,.$ The submanifolds given by
$\det\eta_N=0$ correspond to the components of the subdiscriminant $(1,0,0)\,.$
In terms of the polynomial (\ref{canonicaldef}) these correspond to
further degeneracies amongst the its zeroes.

The induced algebra, Euler vector field and unity
vector field are given by:
\begin{eqnarray*}
\partial_u \star \partial_u & = & [u (u^2 - 2 v^2)/64] \,\partial_u -
[v ( u^2 + v^2 )/64] \,\partial_v \,, \\
\partial_u \star \partial_v & = & -
[v ( u^2 + v^2 )/64] \,\partial_u - [3 u v^2/64] \,\partial_v \,, \\
\partial_v \star \partial_v & = & -[3 u v^2/64] \,\partial_u -
[v ( 4 u^2 + v^2)/64] \,\partial_v\,,\\
&&\\
E_N & = & u \,\partial_u + v \,\partial_v \,, \\
e_N & = & [192 u v /\Delta] \,\partial_u -[64 (u^2 +
v^2)/\Delta] \,\partial_v\,,
\end{eqnarray*}
where $\Delta = v (2 u^2 - v^2)^2\,.$

The $T=T^{(3,1)}$-flow for the $A_3$
Frobenius manifold may be calculated and then restricted
onto the discriminant to yield the hydrodynamic system

\begin{eqnarray*}
u_T & = & (3 u^2 - 3 v^2) u_x - 6 u v v_x \,, \\
v_T & = & -6 u v u_x - ( 3 u^2 + 3 v^2) v_x\,.
\end{eqnarray*}
This may easily be put into Hamiltonian form using the induced
intersection form (up to some overall constants)
\[
\left(
\begin{array}{c}
  u \\
  v
\end{array}
\right)_T =
\left(
\begin{array}{cc}
  1 & 0 \\
  0 & 1
\end{array}\right)\, \frac{d~}{dX}
\,\left(
\begin{array}{c}
  \partial_u h \\
  \partial_v h
\end{array}\right)
\]
where $h=u^4-6 u^2 v^2 - v^4\,.$ This $h$ belongs to a family
given by hypergeometric functions
\[
h_{(1,r)}(u,v) = u^r \,\, {}_2F_{1} \left(
-\frac{r}{2}\,,\frac{1-r}{2}\,,\frac{3-r}{4}\,;\frac{1}{2}
\left(\frac{v}{u}\right)^2\right)\,,
\]
the second family $h_{(2,r)}$ coming the second linearly independent solution of
the corresponding hypergeometric equation.
}
\end{example}

\begin{example}
{\rm
Consider the Frobenius manifold defined by the prepotential and Euler vector field \cite{K,S1}
\begin{eqnarray*}
F & = & t^1 t^2 t^3 - 1/2 t^2 (t^3)^2 + 1/2 t^4 (t^1-t^3)^2 - t^3 e^{t^2} +
e^{t^4} (1+t^3 e^{-t^2}) + 1/2 (t^3)^2 \left(\log t^3 - 3/2\right)\,,\\
E & = & t^1 \partial_1 + \partial_2 + t^3 \partial3 + 2 \partial_4\,.
\end{eqnarray*}
The submanifold corresponding to the limit $t^4\rightarrow-\infty$ is a caustics, since the
polynomial (\ref{canonicaldef}) has a repeated root. However the inverse metric on this
caustic is degenerate:
\[
{}^{(1)}\!g^{ij} = \left(\begin{matrix} 0 & 1 & 0 \\ 1 & 0 & 1 \\ 0 & 1 & 0\end{matrix}
\right)\,.
\]
Thus the ideas in this paper cannot be directly applied. However, the algebra on this
submanifold is still associative, with product
\[
\frac{\partial~}{\partial t^i} \circ \frac{\partial~}{\partial t^j} =
\sum_{k=1}^3 \left.c_{ij}^{~~k} \right|_{t^4 \rightarrow -\infty} \frac{\partial~}{\partial t^k}
\]
so
\begin{eqnarray*}
\partial_i \circ \partial_1 & = & \partial_i \,, \\
\partial_2 \circ \partial_2 & = & - e^{t^2} \partial_2 - t^3 e^{t^2} (\partial_1+\partial_3)\,,\\
\partial_2 \circ \partial_3 & = & - e^{t^2} (\partial_1+\partial_3)\,,\\
\partial_3 \circ \partial_3 & = & - \partial_3 + 1/t^3 \partial_2\,.
\end{eqnarray*}
One may also reduce the corresponding hydrodynamics systems onto this submanifold, but the
systems, while bi-Hamiltonian, have non-trivial Casimirs. Further investigation
of such degenerate caustics, similar to the classical limit of quantum cohomology, requires
further study.
}
\end{example}

\section{Comments}

There are clearly many questions that may be addressed on the structure of submanifolds
in general and of submanifolds in particular. Some of the most interesting concern the connection
with $\tau$-functions are isomonodromy. For an arbitrary integrable system with conserved
densities $h_k$ one may defined the $1$-form
\[
\omega = \sum_i h_i \,dT_{i+1}
\]
which is closed. This then implies the existence of a so-called $\tau$ function
\[
h_k = \frac{\partial\log\tau}{\partial X\partial T_{k+1}}\,.
\]
This types of definition predating more sophisticated definitions based on Grassmannians and
loop groups. On a submanifold $N\subset M$ one may pull back the form, which remains closed,
and hence one may define a $\tau$-functions for the submanifold. For the dispersionless
integrable systems associated with Frobenius manifolds the central object is the
isomonodromic $\tau$-function, denoted $\tau_I\,.$
It would be of interest to see how such an object, and the whole theory of isomonodromy,
behaves on a natural submanifold. One problem is that most of the objects are defined on
$ \mathbb{C}^m \backslash {\rm caustics\,,}$ so various limiting arguments will have to be used to
understand the behaviour of the objects on the caustics themselves. The fact that
Frobenius submanifolds lie in such caustics suggests that this may be possible, at least
in some cases. The singular nature of the $\tau_I$ function on natural submanifolds is
also reminiscent of the work of \cite{AM}, where singularities in $\tau$ functions are
labelled by Young tableaux. This suggest that a general study of the zero/singular set of
$\tau$-functions would be of interest. Intimately connected with the $\tau_I$-function is
the whole question of how one may deform these dispersionless hierarchies \cite{DZ}.

\bigskip

Other interesting questions include:

\begin{itemize}
\item{} To what extent, if at all, does a natural submanifold of a Frobenius manifold define
a (topological/cohomological) quantum field theory?
\end{itemize}
The fact that one has a Frobenius algebra on each tangent space indicates that such
an interpretation may be possible. Also much of the discussion on cohomological field
field theories in \cite{M} may be reproduced without the flatness and condition and the
existence of a prepotential. A simpler question would be to understand how the field
theory corresponding to a Frobenius submanifold is embedded within the larger field theory.

\begin{itemize}
\item{}Do natural submanifolds carry information relevant to ennumerative geometry and quantum
cohomology?
\end{itemize}
In the case of Frobenius submanifolds the information they can contain includes
certain contracted Gromov-Witten invariants of the ambient manifold \cite{M,S2}.
\bigskip

\noindent{\bf Acknowledgements:} I would like to thank Claus Hertling and J\"urgen Berndt for
various useful conversations.


\bigskip

\begin{thebibliography}{9999}

\bibitem[A]{A} {Arnold, V.I., {\sl Singularities of Caustics and Wave Fronts}, Kluwer (1990)}

\bibitem[AF]{AF} {Antonowicz, M and Fordy, A.P., {Coupled KdV equations with multi-Hamiltonian
structures,} Phys. D. {\bf 28} (1987) 345-357.}

\bibitem[AM]{AM} {Adler, M. and van Moerbeke, P., {\sl Birkhoff Strata, B\"acklund Transformations, and
Regularization of Isospectral Operators}, Advances in Math. {\bf 108} (1994) 140-204.}

\bibitem[D]{D} {Dubrovin, B., {\sl Geometry of 2D topological field theories} in {\sl Integrable
Systems and Quantum Groups}, ed. Francaviglia, M. and Greco, S.. Springer lecture
notes in mathematics, {\bf 1620}, 120-348.}

\bibitem[DZ1]{DZ1} {Dubrovin, B. and Zhang, Y., {Extended affine Weyl groups and Frobenius manifolds},
 Compositio Math. {\bf 111} (1998) 167-219.}

\bibitem[DZ2]{DZ} {Dubrovin, B. and Zhang, Y.,
{\sl Bihamiltonian hierarchies in the 2D Topological
Field Theory at One-Loop Approximation}, Commun.Math.Phys. {\bf 198} (1998) 311-361,
{\sl Frobenius Manifolds and Virasoro Constraints}, Selecta Math,. New ser.5 (1999) 423-466,
{\sl Normal forms of hierarchies of integrable PDEs, Frobenius manifolds
and Gromov-Witten invariants}, math/0108160.}

\bibitem[F]{F} {Ferapontov, E.V., {\sl Nonlocal Hamiltonian operators of hydrodynamic type:
differential geometry and applications} in {\sl Topic in topology and mathematical physics},
ed. Novikov, S.P.. Amer. Math. Soc. Transl. Ser. 2 {\bf 170} (1995).}

\bibitem[FP]{FP} {Ferapontov, E.V. and Paplov, M.V., {Quasiclassical limits of couples KdV equations:
Riemann invariants}, Phys. D. {\bf 52} (1991) 211-219.}

\bibitem[FS]{FS} {Fairlie, D.B., and Strachan, I.A.B., {\sl The algebraic and Hamiltonian
structure of the dispersionless Benney and Toda hierarchies}, {Inverse Problems} {\bf 12} (1996) 885-908.}

\bibitem[H]{He} {Hertling, C., {\sl Multiplication on the tangent bundle}, math/9910116.}

%\bibitem[Hi]{Hi} {Hitchin, N.J., {\sl Frobenius manifolds (Notes by David Calderbank)} in
%{\sl Gauge theory and symplectic geometry}, ed. Hurtubise, J. et al., Kluwer Academic
%Publishers (1997). }

\bibitem[HM]{HM} {Hertling, C. and Manin, Yu., {\sl Weak Frobenius manifolds},
Int. Math. Res. Notices {\bf 6} (1999) 277-286.}

\bibitem[K]{K} {Kodama, Y., {\sl Dispersionless integrable systems and their solutions}, in
{\sl Integrability: the Seiberg-Witten and Whitham equations}, ed. Braden, H.W. and Krichever, I.M.,
Gordon and Breach (2000).}

\bibitem[M]{M} {Manin, Y.I., {\sl Frobenius manifolds, Quantum
Cohomology, and Moduli Spaces}, A.M.S. Colloquium, Vol.47 (1999).}

\bibitem[Mo]{Mo} {Mokhov, O.I., {\sl Compatible and Almost Compatible Pseudo-Riemannian Metrics},
math.DG/0005051 }

\bibitem[Sh]{Sh} {Shcherbak, O.P., {\sl Wavefront and reflection
groups},
Russ. Math. Surv. {\bf 43:3} (1988) 149-194.}

\bibitem[St1]{S1} {Strachan, I.A.B., {\sl Degenerate Frobenius manifolds and the
bi-Hamiltonian structure of rational Lax equations}, {J. Math. Phys.} {\bf 40} (1999) 5058-5079.}

\bibitem[St2]{S2} {Strachan, I.A.B., {\sl Frobenius submanifolds},
{J. Geom. Phys.}, {\bf 38} (2001) 285-307.}

\bibitem[St3]{S3} {Strachan, I.A.B., {\sl Frobenius manifolds and bi-Hamiltonian
structures on discriminant hypersurfaces}, to appear in the proceedings of the
conference {\sl Integrable Systems in Differential Geometry} held in Tokyo, July 2000.}
%{\sl Integrable Systems, Topology
%and Physics}, A.M.S. Contemporary Mathematics series. }

\bibitem[Z]{Z} {Zuber, J.-B., {\sl On Dubrovin topological field theories}, Mod. Phys. Lett.
{\bf A9} (1994) 749-760.}

\end{thebibliography}
\end{document}